\documentclass[12pt,a4paper]{article}

\usepackage{amsmath,amssymb}
\usepackage{color}
\usepackage{hyperref}
\usepackage{graphicx}
\newcommand{\R}{\mathbb{R}}

\newcommand{\Z}{\mathbb{Z}}

\newcommand{\cqfd}{{\nobreak\hfil\penalty50\hskip2em\hbox{}\nobreak\hfil
$\square$\qquad\parfillskip=0pt\finalhyphendemerits=0\par\medskip}}

\newcommand{\tcb}{\textcolor{blue}}

\newcommand{\D}{\Delta}

\newcommand{\n}{\nabla}

\newcommand{\p}{\partial}

\newtheorem{prop}{Proposition}
\newtheorem{defi}{Definition}

\newcommand{\e}{\epsilon}
\newcommand{\va}{\varphi}

\newtheorem{theorem}{Theorem}

\newtheorem{remarka}{Remark}

\newtheorem{lemme}{Lemma}

\title{  Fractional $BV$ solutions       
for  $2\times 2$  systems of  conservation laws 
 with  a  linearly degenerate field}
  \author{Boris Haspot  \thanks{Universit\'e Paris Dauphine, PSL Research University, Ceremade, Umr Cnrs 7534, Place du Mar\' echal De Lattre De Tassigny 75775 Paris cedex 16 (France), haspot@ceremade.dauphine.fr  } 
,  St\'ephane  Junca \thanks{Universit\'e C\^ote d'Azur, Inria \& CNRS, LJAD,  France, Stephane.JUNCA@univ-cotedazur.fr}
}
\date{\empty}
\begin{document}

\maketitle
\begin{abstract}
The class of $2\times 2$  nonlinear hyperbolic systems with one genuinely nonlinear  field 
and one linearly degenerate field are considered.
Existence of global weak solutions for small initial data in fractional BV spaces $BV^s$ is proved. The exponent $s$ is related to the usual fractional Sobolev derivative. 
Riemann invariants $w$ and $z$  corresponding respectively to the genuinely nonlinear component and to  the linearly degenerate component play different  key roles in this work. 
We obtain  the existence of  a global weak solution provided  that  the initial data  written in Riemann coordinates 
$ (w_0,z_0)$ are small in $ BV^s \times L^\infty $, $1/3 \leq s<1$.
The restriction on the exponent $s$ is due to a  fundamental result of P.D. Lax, the variation of the Riemann invariant $z$ on the Lax shock curve depends in a cubic way of the variation of the other Riemann invariant $w$.  
\end{abstract}

\tableofcontents
\section{Introduction}

In this paper, we study general $2\times 2$  hyperbolic systems of the form:
\begin{equation}\begin{cases}
\begin{aligned}
&\p_t U+\p_x F(U)=0\\
&U(0,\cdot)=U_0
\end{aligned}
\end{cases}
\label{1}
\end{equation}
with $U(t,x)\in \Omega \subset \R^2$ an open set, $(t,x)\in\R^+\times\R$. $F$ is the flux of the system and it is regular from $\R^2$ to $\R^2$. We assume that the system is strictly hyperbolic  {on $\Omega$}, it means that $DF(U)$ has two different eigenvalues $\lambda_1$ and $\lambda_2$. Without any restriction we can assume that $\lambda_1<0<\lambda_2$,  {reducing if necessary the open set $\Omega$}. It implies in particular that we have a basis of eigenvector of unit norm $(r_1(U),r_2(U))$ for any $U\in \Omega \subset \R^2$.
In the sequel we will only be interested in the case of a 1 genuinely nonlinear field and a 2 linearly degenerate field. In particular, it means that for every $U\in  \Omega \subset \R^2$ we have 
\begin{equation}
\n\lambda_1(U)\cdot r_1(U)\ne 0\;\;\mbox{and}\;\;\n\lambda_2(U)\cdot r_2(U)= 0.
\label{conf}
\end{equation}
\paragraph{Examples} We wish now to give some examples of strictly hyperbolic system 
with a  genuinely nonliear field and a linearly degenerate field.
\begin{itemize}
\item  The classical chromatography system  \cite{Br,Da} when the velocity is known which is the case for the liquid chromatography.  
\item     The Keyfitz-Krantzer system \cite{KK} has this structure, it is maybe the 
first and the most famous known. It is related to some problem of
nonlinear elasticity.
\item 
The $2\times 2$ Baiti-Jenssen system \cite{BaJe97} with a genuinely nonlinear field. The Baiti-Jenssen systems  arise in models for porous media, traffic and
gas flows. 
\item  The  Aw-Rascle  system is well known for traffic flow \cite{AR}.
\\

The  four first examples are Temple systems \cite{Serre,Temple83, Temple}. Such systems satisfy a maximum principle which is not generally true for systems of conservation laws. 
Now the following list provides examples that are not  Temple systems.
\\

\item  
Colombo and Corli consider the class of $2\times 2$ system with genuinely nonlinear field and a Temple field \cite{CoCo}.  They prove existence of solutions for large $BV$ data associated to the Temple component.      A linear degenerate field is an example of Temple field, the rarefaction and shock curve coincide \cite{Br}. Such assumption is not enough to have a Temple system. One interest of our work is to prove existence in $BV^s$, so, with possible infinite total variation in $BV$.

\item The chromatography system with a sorption effect  \cite{BGJ2}     
is   a chromatography system with  a non constant  and unknown velocity. 
   This   system is generally not a Temple system    \cite{BGJ4}. 
 \item  We mention also some triangular systems  with a  transport equation \cite{Triang}. This class of systems generalizes the previous one when it is written in 
Lagrangian coordinates \cite{BGJP}.
\end{itemize}
In this paper we would like to extend, for  $2\times 2$  systems with a genuinely nonlinear filed and a linearly degenerate one, the famous result of Glimm \cite{Glimm} concerning the existence of global weak solution for the strictly hyperbolic system with small initial data $u_0$ in $BV$. Indeed we would like to enlarge the set of initial data by working with $u_0$ belonging  to $BV^s$ with $0<s<1$, $BV^1=BV$. We are now going to give a definition of the fractional BV spaces called $BV^s$. We refer also to Bruneau \cite{Bru} for more details.
\begin{defi}[$TV^s$ variation] \label{defcle}
We say that a function $u$ is in $BV^s(\R)$ with $0<s\leq 1$  and $p=1/s \geq 1$ if
$TV^su<+\infty$ with:
\begin{equation} \label{def:TVs}TV^s u:=\sup_{n \in \mathbb{N},\;  x_1 < \cdots< x_n}\sum_{i=1}^n |u(x_{i+1})-u(x_i)|^{p}
\end{equation}
The associated semi-norm of  the $TV^s$ variation is, 
 \begin{equation} \label{def:semi-normBVs}
 |u|_{BV^s}:=   (TV^s u)^s
 \end{equation}
and a norm is 
 \begin{equation} \label{def:normBVs}
 \|u\|_{BV^s}:= \|u\|_{L^\infty}+     |u|_{BV^s}
 \end{equation}
\end{defi}
In the same way, $TV^su(I)$ is defined as the $TV^s$ variation of the function $u$ on the set $I$. 
We note that it is clear that for any $s\in]0,1]$, $BV^s(\R) \subset L^\infty(\R)$  \cite{MO}.  Moreover, if $u$ belongs to $L^1(\R)$ then the semi-norm $BV^s$ is a norm.  This is due to the fact that a $BV^s$ function has limts at $\pm  \infty$ and, for a $L^1(\R)$ function, these limits are $0$.
  For $0<s_1<s_2\leq 1$, we also  have $BV^{s_2}\subset BV^{s_1}$  \cite{BGJ6}.
  The $TV^s$ variation was called the $p-$variation with $p=1/s$ in \cite{MO}.  We prefer to use the notation $TV^s$  since it is related to the Sobolev exponent ``s''. Indeed, $BV^s_{loc}$ is close from $W^{s,p}_{loc}$ but remains different  \cite{BGJ6}, indeed the $BV^s$ functions are regulated functions  \cite{MO} as $BV$ functions.

\begin{prop}[$BV^s$ functions are regulated functions \cite{MO}]
If $u\in BV^s$ with $0<s\leq 1$ then $u$ admits only a countable set of discontinuity. Futhermore for every $x\in\R$, $u$ admits a limit on the left and on the right in $x$.
\label{propimp}
\end{prop}
We would like now to motivate the use of the $BV^s$ spaces for the study of hyperbolic systems. Actually the most of the results on the existence of global weak solution for $2 \times 2$ hyperbolic systems concerned the $L^\infty$ space and the $BV$ space. In order to tackle this problem, there exists essentially two different approaches, the first one was developed by Glimm in the 60s \cite{Glimm}. He proved for a general $n \times n$ strictly hyperbolic system with genuinely non linear field or linearly degenerate field the existence of global weak entropy solution provided that the initial data is small in $BV$. The main difficulty of the proof consists in controlling the $BV$ norm of the solution all along the time, indeed Glimm has observed that the $BV$ norm can increase after each interaction between the nonlinear waves. In order to estimate this gain in $BV$ norm after each interaction, Glimm has introduced a quadratic functional 
 which described the  interactions between the nonlinear waves and which allows to evaluate the $BV$ norm of the solution all along the time.
This result has been extended in the 90s by Bressan and his collaborators  \cite{Br,BrGoat,BreLe} where they proved the uniqueness of Glimm solution (provided that $U_0$ belongs also to $L^1(\R)$) in a suitable class of solution which takes into account in particular the Lax conditions for the shocks. The main ingredient to do this is to prove that the wave front tracking algorithm (we refer to \cite{Br} for the definition of the wave front tracking for general $n \times n$ systems ) generates a Lipschitz semigroup in $L^1$ \cite{Br,BrCo}. We recall in particular that the solutions which are issue of the wave front tracking method and which are determined via a compactness argument are the same as the solution coming from the Glimm scheme  \cite{Br}.\\
The second approach was initiated by Di Perna  \cite{Dip1,Dip2} at the beginning of the 80s using the so called compensated compactness which was introduced by Tartar  \cite{Tar}. Roughly speaking this method can be applied for $2 \times 2$ strictly hyperbolic systems with two genuinely non linear fields (see also Serre \cite{Serre}) when the initial data $U_0$ is assumed to belong to $L^\infty(\R)$. The case of the isentropic Euler system has been particularly studied and we refer to \cite{Chen,L1,L2,LF}. We observe then that this method allows to deal with more general initial data as $U_0\in BV$, however there is generally no result of uniqueness for these solutions. In particular it seems complicated to select the solution via the Lax conditions on the shocks since we can not give any sense of traces along a shock for such solutions since they belong only to $L^\infty_{t,x}$.\\
In the $2 \times 2$ case, Glimm in \cite{Glimm} has obtained a better result of existence of global weak solution
inasmuch as he can deal with large initial data $U_0$ in $BV$ provided that the $L^\infty$ norm of
$U_0$ is sufficiently small. It is due to the fact that after an interaction between waves the variation of the $BV$ norm  has  a cubic  {order} in terms of the incoming strengths of the waves which interact (in the general case $n \geq 3$,  this order is only quadratic). This result is a consequence of the existence of Riemann invariants for $2 \times 2$ systems. Later on, this result has been extended by Glimm and Lax in \cite{GL} to the case of small $L^\infty$ initial data when the fields are genuinely non linear. We refer also to the recent work of Bianchini, Colombo and Monti  \cite{BCM}. To do this, the authors use new Glimm functionals to control the $L^\infty$ norm combined with the method of backward characteristics. In addition they proved a new Oleinik inequality (which is generally restricted to the scalar conservation law with genuinely non linear flux) for this $2 \times 2$ system which gives sufficient compactness informations in order to pass to the limit in the wave front tracking. 
 \\ \\
{$BV^s$ spaces are intermediary spaces between $L^\infty$ and $BV$ spaces, see \cite{BGJ6} or  the definition \ref{defcle}. We note in particular that the $BV^s$ spaces admits  functions with shocks, from this point of view these spaces are suitable for dealing with hyperbolic systems of conservation laws. Indeed it is well known that the solution of an hyperbolic system can admits shock in finite time even if the initial data is arbitrary regular. In addition (see  \cite{MO}  and   the proposition \ref{propimp} below),   $BV^s$ functions admit a notion of ``traces'' as for $BV$  functions (this is of course not the case for $L^\infty$  functions).}
This notion of trace is essential in the result of uniqueness of Bressan et al \cite{Br,BrGoat,BreLe}. Indeed, it gives a sense to the notion of the Lax entropy criterion which enables to select a unique solution (in the results of Bressan \& al. 
a tame oscillation  condition is also required). It would be then interesting to prove the existence and the uniqueness of global weak solution for initial data $U_0$ in $BV^s$ with $0<s<1$ for strictly hyperbolic systems. It would improve in particular the results of existence  of Glimm \cite{Glimm} inasmuch as the initial data $U_0$ would be less regular as $BV$. In addition, working with $BV^s$ gives a chance to extend the results of uniqueness of Bressan \& al.   \cite{Br,BrGoat,BreLe} since the notion of trace remains relevant.
\\
\\
In this paper, we will only focus our attention on the existence of global weak solutions for small initial data in $BV^s$. Before giving and describing our main results,
we would like now to recall some results using the $BV^s$ space in the framework of conservation laws.  {For scalar conservation law, the entropy solution corresponding to an initial data belonging in $BV^s$ remains in $BV^s$ for all time \cite{BGJ6}. Moreover, this result is sharp \cite{allt}. }
 We would like to point out that the $BV^s$ space is also naturally used to describe the regularizing effects of scalar conservation laws. From \cite{K}, it is known that there exists unique global solution for scalar conservation laws when $U_0$ belongs to $L^\infty$. The most famous regularizing effect concerns the uniformly convex flux where the solution becomes instantaneously $BV_{loc}$, this is a direct consequence of the so called Oleinik inequality. For a convex flux with a power law degeneracy, the authors in $\cite{BGJ6}$ show an optimal regularizing effects on the solution $u$ inasmuch as the solution $u$ becomes instantaneously $BV^s_{loc}$ with $s$ depending on the power law of the flux. { For  a general nonlinear convex flux locally Lipschitz,   the solution belongs for positive time in a generalized  $BV$ space related to the nonlinearity of the flux, see \cite{Lipconv}.
 These results have been extended  for a nonlinear non convex flux, at least $C^3$,   in a generalized  $BV$ space, and for a more regular flux  with polynomial degeneracy  in the optimal $BV^s$ space by Marconi in \cite{EM,EMproc}.}\\
 \\
The most of the result dealing with $BV^s$ initial data concerned scalar conservation laws. Indeed it is a priori delicate to prove the stability of the $BV^s$ norm all along the time, the $BV^s$ norm is indeed more complicated to compute than the $BV$ norm. Indeed when we apply a Glimm scheme, in order to know the $BV$ norm after an interaction between waves, it is sufficient to estimate locally the strength of the new outgoing waves since we recover the complete $BV$ norm by summing the absolute value of the different strength on  all the euclidean space. In particular using the triangular inequality, we do not need to select subdivisions of the euclidean space in a accurate way in order to control the $BV$ norm. It is not the case for the $BV^s$ norm
which  is reached for particular optimal subdivisions. It implies in particular that for computing the $BV^s$ norm after a waves interaction, it is not sufficient to knows only the values of the outcoming strength.
In the scalar case, the analysis is simpler since after each interactions, there exists some zone of monotonicity for the Riemann problem making the analysis simpler to compute the $BV^s$ norm \cite{BGJ6}.\\

 In this paper, we would like to extend the analysis of \cite{BGJ6} to the case of $2 \times 2$ strictly hyperbolic system
with one genuinely nonlinear field and one linearly degenerate field which corresponds to the case described in (\ref{conf}). This case is a particular case of the theory of Glimm \cite{Glimm} on the existence of global weak solution for initial data $U_0$ in $(BV(\R))^2$ with a large $BV$ norm provided that the $L^\infty$ norm is sufficiently small. We would like also to mention that others authors  have yet obtained existence of weak entropy solution  for small $L^\infty$ data and large $BV$ norm when an eigenvalue is linearly degenerate \cite{BCM,Da} or a Temple eigenvalue \cite{CoCo}. 
We extend the results of Glimm by proving the existence of global weak solution for small initial data with $(w_0,z_0)$ belonging in $BV^s\times L^\infty$ with $\frac{1}{3}\leq s<1$. Here $(w,z)$ are the solution of the system (\ref{1}) that we consider in Riemann coordinates respectively in terms of the 1 genuinely nonlinear field and the 2 degenerate field. To do this,  we follow the classical method which consists in introducing a wave front tracking with $(w^\nu,z^\nu)$ corresponding to the approximate solutions $(U^\nu)$ of the wave front tracking  written in Riemann coordinates and $\nu\rightarrow +\infty$ the parameter associated to the wave front tracking.  We are then interested in proving that $(U^\nu)_{\nu>0}$ converges to $U$ a solution of the system (\ref{1}).
The main difficulty consists in proving uniform $BV^s$ estimates on $(w^\nu,z^\nu)$ and next in verifying that the wave front tracking is well defined for any time $t>0$. The end of the proof requires to establish compactness argument in order to verify that $U^\nu$ converges to $U$ a solution of the system (\ref{1}) (here $U^\nu$ is the approximated solution associated to the wave front tracking written in physical coordinates and not in Riemann coordinates).\\ 
More precisely we show that the $BV^s$ norm of $w^\nu$ is uniformly conserved all along the time essentially because the {waves interactions} do not increase the $BV^s$ norm for $w^\nu$.
The proof is reminiscent of the scalar case for a convex flux. However it is more complicated to control uniformly the $L^\infty$ norm of $z^\nu$. Indeed the norm of $z^\nu$ can increase after two types of interactions, interaction between 1-shocks and interaction between 1-shock and a 2-contact discontinuity.
To do this, we consider the $L^\infty$ norm of $z^\nu$ along  any  forward generalized 2-characteristic and we observe that this $L^\infty$ norm depends on the $BV^{\frac{1}{3}}$ norm of $w_0$. Indeed the $L^\infty$ norm of $z^\nu$ along a  forward generalized 2-characteristic increases only when the  forward generalized 2-characteristic meets a 1-shock, furthermore this increase depends on the cubic strength in { the variation} of $w^\nu$ on this 1-shock ( it is important to point out that this increase is directly related to the regularity of the Lax shock curve).
It suffices then to follows these 1-shock in a backward manner in order to estimate the $L^\infty$ norm of $z^\nu$ in terms of $\|w_0\|_{BV^{\frac{1}{3}}}. $ It explains why we need to assume that $w_0$ is in $BV^s$ with $s$ at least equal to $\frac{1}{3}$. The last step of the proof consist in proving that $z^\nu$ converges to $z$ up to a subsequence in $L^1_{loc,t,x}$. This part is a priori delicate since we have only a uniform control of the $L^\infty$ norm of $z^\nu$. We observe however that we have additional regularity property if we study the unknown $z^\nu_L(t,x)=z^\nu(t,\gamma^\nu_2(t,x))$ with $\gamma^\nu_2(t,x)$ the  forward generalized 2-characteristic such that $\gamma_2^\nu(0,x)=x$. Here $z^\nu_L$ is the value of $z$ in Lagrangian coordinates, following the same idea as for the control of the $L^\infty$ norm of $z^\nu$, we can prove that $z_L^\nu$ is uniformly bounded in $L^\infty_x(BV_t)$ and that the speed of propagation is finite. It is then sufficient to prove that $z^\nu_L$ converges up to a subsequence strongly in $L^1_{loc,t,x}$, we prove next that the Lagrangian transformation $(t,\gamma_2^\nu(t,x))$ is a uniformly bi-Lipschitz homeomorphism in $\nu$ what is sufficient to ensure that $z^\nu$ converges also strongly in $L^1_{loc,t,x}$. It allows to conclude that the solution $(U^\nu)$ converges to a solution of (\ref{1}). 
\section{Presentation of the results}
We would like to state now our main result. 
For this purpose we use  a  distinguished coordinate
system called Riemann invariants, which in general exists only for $2\times 2$ systems (chapter 20, \cite{Smoller}). This allow to perform a nonlinear diagonalization of the hyperbolic system for smooth solutions. This diagonalization is not valid for discontinuous solutions but the Riemann invariants have still some advantages.  The Riemann  problem and the interaction of waves is also simpler to study in these coordinates than in the initial coordinates.   The following notations $(w,z)$ are chosen for the Riemann invariants.
Thanks to Riemann \cite{Smoller}, there exists a change of coordinates $U \mapsto (w,z)=(w(U),z(U))$ (here $(w,z)=(w(U),z(U))$ is a standard abuse of notations), reducing the open set $\Omega$  if necessary, such that
\begin{align}
\nabla  w \cdot r_2 = 0,   & &  \nabla  z \cdot r_1 = 0.
\end{align}
In all the sequel, $U$ is written in this system of coordinates.  In particular, the initial data $U_0$ of the system \ref{1} reads $w_0=w(U_0)$ and $z_0=z(U_0)$.
Our main theorem states as follows.
\begin{theorem} [Existence in $BV^{1/3} \times L^\infty $]
\label{thm:1/3,0}
Let $w_0\in BV^s(\R)$ with $\frac{1}{3} \leq s \leq 1$  and $z_0\in L^\infty(\R)$ then there exists  $\e_0>0$ such that if:
$$\|w_0\|_{BV^s}+\|z_0\|_{L^\infty}\leq\e_0$$
then there exists a global weak solution $U$ for the system (\ref{1}). The Riemann coordinates $(w,z)$ belong to $L^\infty(\R^+,BV^s\times L^\infty)$. 
\\
Moreover, the Riemann invariant $z$ can be decomposed  in the Lagrangian coordinates associated to the linearly degnerate field as follows: 
\begin{align}\label{zdecomposition}
z(t,\gamma_2(t,x)) & =z_0(x) +\eta(t,x)
\end{align}
where $\lambda_2$ is the linearly degenerate eigenvalue which depends only on $w$ and $\gamma_2$ represents  the generalized 2-characterisics, 
\begin{align*}
\begin{cases}
\displaystyle
\frac{d \gamma_2}{d\, t} (t,x) & = \lambda_2(w(t,\gamma_2(t,x))) \\
\gamma_2(0,x) &= x  
\end{cases}
\end{align*}
with $\eta    \in L^\infty_x(\R, BV_t(\R^+)) \cap Lip_x(\R, L^1_{loc,t}(\R^+))$.
\label{theo2}
\end{theorem}
\begin{remarka}
In this Theorem, we assume that $(w_0,z_0)\in BV^s\times L^\infty$, $\frac{1}{3}\leq s<1$.  A sufficient condition  on the initial physical coordinates $U_0$ to ensure such regularity is to take  $U_0 \in  BV^s$. 
\end{remarka}
\begin{remarka}
The decomposition of $z$ in \eqref{zdecomposition} provides also a stability results in $BV^\sigma$ for all $0 \leq \sigma \leq 1$. 
This means that $z$ is in $L^\infty(\R^+,BV^\sigma)$ if $z_0 \in BV^\sigma$.
\end{remarka}
This decomposition \eqref{zdecomposition} has already been obtained  as a factorization of the gas velocity for a chromatography system  (Theorem 7.2 in \cite{BGJ4}). \\
\\
Up to our knowledge, Theorem \ref{thm:1/3,0} is the first general result concerning the stability of the $BV^s$ norm in the framework of strictly hyperbolic systems, except for some particular physical systems \cite{BGJP,JuLo5}. It extends in particular the results of Glimm \cite{Glimm} since the initial data are not necessary $BV$. Furthermore if we compare this result with the works of Glimm, Lax and Bianchini,Colombo, Monti (see \cite{GL} and \cite{BCM}) which deal with initial data in $L^\infty$, we can only say that the framework is different. Indeed in \cite{GL,BCM} the two fields are genuinely nonlinear, in particular the authors extend the Oleinik inequality to their case what allows them to
obtain sufficient information in terms of compactness to pass to the limit respectively in their scheme and their wave front tracking. In our case, the technics are quite different especially on the arguments of compactness which enables us to consider the limit of the approximated solutions $(w^\nu,z^\nu)$ which are issue of the wave front tracking. Indeed, an observation is to remark that the solution $z^\nu$ can be splited into $z_0$ the initial data and a function $\eta^\eta$ which is more regular as $z^\nu$ itself. We can then pass to the limit in $\nu$ for $\eta^\nu$ in the wave front tracking.
\begin{remarka}
 For $2\times 2$ Temple system with a genuinely nonlinear field and a linearly degenerate field, these results improve the classical existence in $BV$. The existence for $L^\infty$ data in \cite{BrGoat} needs that all fields are genuinely nonlinear. 
\end{remarka}
\begin{remarka}
Since the $BV^s$ norm has a notion of trace it would be interesting to prove the uniqueness of the solution.
\end{remarka}
In the section \ref{secinitial}, we start by giving one useful Lemma to compute in a simple way the $BV^s$ norm of a sequence $(u_n)_{n\in\mathbb{N}}$.
In the section \ref{sec:WFT}, we describe the Lax curve and the different interactions between 1 and 2-waves. Furthermore we define a simplified wave front tracking  well adapted to our case which concerns a one genuinely nonlinear field and one linearly degenerate. In the section \ref{sec4}
we prove the Theorem \ref{thm:1/3,0}.
\section{Local monotonicity and computation of $TV^s u$}
\label{secinitial}
The computation of the $TV^s$ variation can be more complicate than the usual $TV$ variation \cite{BGJ6}. Here, we prove in a self contained and  in a detailed way an important tool  to  compute the $BV^s$ norm for a function which is piecewise constant. 
\\
\\
 For a sequence  $(u_n)_{n=1,\cdots,N}$ , a subdivision $\sigma$  is considered as a subset of  $\{1,\cdots,N\}$ or  as an increasing application from  $\{1,\cdots,|\sigma|\}$ to $\{1,\cdots,N\}$ where $|\sigma|\leq N$ is the the cardinal of the subdivision.  This means that the subdivision $\sigma$ can be written in terms of  the bijection $\sigma$ as $\{\sigma(1),\cdots,\sigma(|\sigma|) \}$
\begin{lemme}
\label{lemcru}
Assume that we have a sequence $(u_n)_{n=1,\cdots,N}$ with $i_0\in\{2,\cdots,N-1\}$:
$$u_{i_0-1}< u_{i_0}< u_{i_0+1}$$
and assume that for $0<s\leq 1$:
$$TV^su=\sum_{ 1 \leq i  <|\sigma|} |u_{ \sigma(i+1)}-u_{\sigma(i)}|^{\frac{1}{s}}$$
with $\sigma$ a subdivision of $\{1,\cdots,N\}$ which is optimal to compute  $TV^s u$.
Then $i_0$ is not in the subdivision $\sigma$. In particular if we set $v=({u_1,\cdots,u_{i_0-1},u_{i_0+1},\cdots,u_N})$ then we have:
$$TV^s u=TV^sv.$$
The conclusion is the same if we have:
$$u_{i_0-1}> u_{i_0}> u_{i_0+1}.$$
\end{lemme}
{\bf Proof:} We assume here by absurd that $i_0\in\sigma$ with $\sigma$ an optimal subdivision for the $BV^s$ norm of $(u_n)_{n=1,\cdots N}$. Furthermore the first term of the subdivision $\sigma$ before $i_0$ is $i_{-1}$ and the next one is 
$i_1$, we have in particular $i_{-1}<i_0<i_1$. It is important to note that for example $i_{-1}$ exists, indeed if we assume that $i_0$ is the first term of the subdivision $\sigma$ it is easy to observe that $\widetilde{\sigma}=\{i_0-1\}\cup \sigma $ is a subdivision which produces a larger $BV^s$ norm as the subdivision $\sigma$. The same argument is also true for $i_1$.
We are now going to consider different cases. We note in the sequel $p=\frac{1}{s}$.
\begin{itemize}
\item $i_{-1}<i_{0}-1$ and $i_1>i_0+1$.\\
In this case if we have $u_{i_{-1}}<u_{i_0}<u_{i_1}$ then $i_0$ can not be in the subdivision since:
\begin{equation}|u_{i_0}-u_{i_{-1}}|^p+|u_{i_0}-u_{i_{1}}|^p< |u_{i_{-1}}-u_{i_{1}}|^p.
\label{rest1}
\end{equation}
Indeed the subdivision $\sigma$ would be not optimal for the $BV^s$ norm.\\
If we have $u_{i_{-1}}<u_{i_0}$, $u_{i_1}\leq u_{i_0}$, then we observe that:
$$|u_{i_{-1}}-u_{i_0}|^p+|u_{i_1}-u_{i_0}|^p< |u_{i_{-1}}-u_{i_0+1}|^p+|u_{i_1}-u_{i_0+1}|^p.$$
It implies that $\sigma$ is not optimal for the $TV^s$ norm since the subdivision $\widetilde{\sigma}$  $\cdots, i_{-1},i_0+1,i_1,\cdots$ produces a larger $BV^s$ norm.\\
The other cases $u_{i_{-1}}>u_{i_0}>u_{i_1}$ and $u_{i_{-1}}>u_{i_0}\leq u_{i_1}$ can be treated similarly, in particular for the last case it suffices to consider the subdivision $\widetilde{\sigma}$ defined by $\cdots, i_{-1},i_0-1,i_1,\cdots$. We note that the case $u_{i_1}=u_{i_0}$ has no interest since we can omit in the subdivision $\sigma$ the term $i_0$.
\item $i_{-1}<i_0<i_0+1=i_1$ with $i_{-1}<i_{0}-1$.\\
If $u_{i_{-1}}<u_{i_0}$ then $i_0$ is not in the subdivision again using (\ref{rest1}).\\
 If $u_{i_{-1}}>u_{i_0}$ then we have:
$$|u_{i_{-1}}-u_{i_0}|^p+|u_{i_0+1}-u_{i_0}|^p< |u_{i_{-1}}-u_{i_0-1}|^p+|u_{i_0-1}-u_{i_0+1}|^p.$$
It means that $\sigma$ is not an optimal subdivision for the $BV^s$ norm, indeed we can consider the subdivision $\widetilde{\sigma}$ $\cdots, i_{-1},i_0-1,i_0+1,\cdots$ which produces a larger $BV^s$ norm.
\item We proceed similarly for the last case  $i_{-1}=i_0-1<i_0$.
\end{itemize}
In conclusion in all the case $i_0$ can not belong to $\sigma$, it concludes the proof.
\cqfd
\section{A wave front tracking algorithm}
\label{sec:WFT}  
In this section, the wave front tracking  algorithm (WFT) used is presented to solve the initial value problem  \eqref{1}. Simplifying the (WFT) is useful to simplify the estimates on the approximate solutions \cite{BaJe98}.
 Taking advantage of the linearly degenerate field, we present a simpler wave front tracking  algorithm (WFT) as the one used  by Bressan and Colombo  in \cite{BrCo}  for general $2\times 2$  systems.   For this purpose, the Riemann problem and the interaction of waves is first studied.  
 \\
 In the sequel, we denote by $A(U)$ the $2 \times 2$ hyperbolic matrix $DF(U)$ and without loss of generality by $\lambda_1<0<\lambda_2$ its eigenvalues and by $l_1$, $l_2$ (respectively $r_1$, $r_2$) its left (repectively right) eigenvectors, normalized so that:
 $$
 \begin{aligned}
 &\|r_i(U)\|=1,\;\langle l_j(U);r_i(U)\rangle=\delta_{i,j},\;i,j=1,2.
 \end{aligned}
 $$
Furthermore we recall  that for all $U \in \Omega$(\ref{conf}):
\begin{equation}
\n\lambda_1(U)\cdot r_1(U)\ne 0\;\;\mbox{and}\;\;\n\lambda_2(U)\cdot r_2(U)= 0.
\label{confa}
\end{equation}
In the sequel we assume that we have $ \n\lambda_1(U)\cdot r_1(U)>0$ for all $U$ in $\Omega$.
Furthermore for $\Omega=B(0,r)$ sufficiently small, we have.
\begin{equation}
\sup_{U\in B(0,r)}\lambda_1(U)<0<\inf_{U\in B(0,r)}\lambda_2(U).
\end{equation}
\subsection{Riemann invariants and Lax curves}

An important feature for $2\times 2$ systems is the existence (at least locally)  of coordinates in the state space, the  Riemann  invariants. 
 All properties of  the solutions  $U$ are stated in the Riemann invariants coordinates.
We call $w$ and $z$ the Riemann invariant associated to genuinely nonlinear (GNL) eigenvalue   $\lambda_1$ and the linearly degenerate one $\lambda_2$. More precisely we have $\n w(U)\cdot r_2(U)=0$ and $\n z(U)\cdot r_1(U)=0$ for any $U\in \Omega \subset \R^2$.    When $U$ is a solution of system \eqref{1}, $(w(t,x),z(t,x)))$ denotes  $(w(U(t,x)),z(U(t,x)))$. With this usual notation, the Riemann invariant satisfies with here by abuse of notation $\lambda_1(w,z)=\lambda_1(U)$ and $\lambda_2(w,z)=\lambda_2(U)$:
$$
\begin{cases}
\begin{aligned}
&\p_t w+\lambda_1(w,z)\p_x w=0\\
&\p_t z+\lambda_2(w)\p_x z=0.
\end{aligned}
\end{cases}
$$
\begin{remarka}
Notice that      $ \partial_w \lambda_1 > 0 $, since the first field is GNL, and  
 $\lambda_2$  is independent of $z$  since the second field is linearly degenerate (see Theorem 8.2.5 \cite{Da}). That is   $ \partial_z \lambda_2 =0 $ and  $ \lambda_2$  depends only on $w$. 
\end{remarka}
The map $U\rightarrow (w(U),z(U))$ is a local diffeomorphism and we can assume that the origin in $U$ coordinates corresponds to the origin in $(w,z)$ coordinates. We are going now to define the Lax curves in these new coordinates $(w,z)$. For a fixed state $U_-$, the Lax curves describe the set of state $U_+$ such that the Riemann problem with the left state $U_-$ and the right state $U_+$ is a simple wave \cite{Smoller}. For each $U_-$ fixed there are two Lax curves, one for the 1-waves and another one for the 2-waves.


\paragraph{Lax Curves}

\begin{figure}
\includegraphics [scale=0.6]{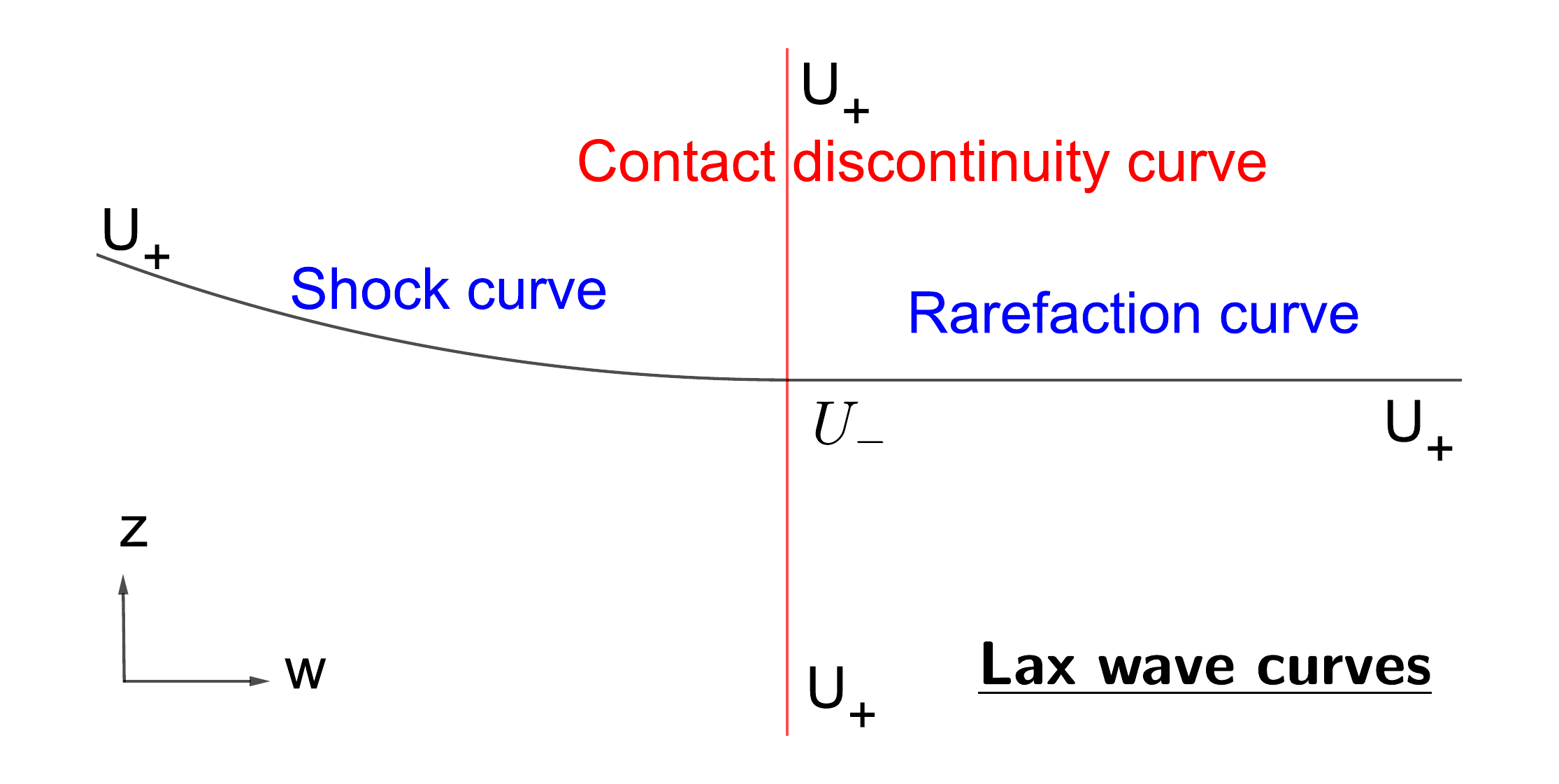}
\caption{Lax waves curves  where the state on the left of the wave $U_-=(w_-,z_-)$ is fixed. $U_+$ is the right state connected by a 1-wave when $w$ varies, $w_+\neq w_-$, or a 2-wave when $w_+=w_-$ is constant. }
\label{fig:Lax-curves}
\end{figure}

The picture (see the figure \ref{fig:Lax-curves}) of the Lax curves $L(U_-)$ will be used systematically throughout the paper. 
The convexity (or the concavity) of the shock curve  is constant  at least on a small  ball $B(0,r)$ with $r>0$ and does not depend on the point $U_-$ and $U_+$ which are both in $B(0,r)$.  
 It simply corresponds to fix a sign of the  third derivative of the  shock curve (we refer to \cite{Da} Theorem 8.2.3).  The other sign can be studied in a similar way.  

\noindent For the rarefaction $R_1$ we have:
\begin{equation}
\begin{cases}
\begin{aligned}
&w=w_-+\sigma\;\;\mbox{with}\;\sigma\geq 0.\\
&z=z_-
\end{aligned}
\end{cases}
\label{aR1}
\end{equation}
For the $S_1$ shock we have:
\begin{equation}
\begin{cases}
\begin{aligned}
&w=w_-+\sigma\;\;\mbox{with}\;\sigma\leq 0.\\
&z=z_-+ \mathcal{O}(\sigma^3)
\end{aligned}
\end{cases}
\label{aS1}
\end{equation}
Notice, with the choice of the convexity for the shock curve, $z$ increases  through a shock wave.  (For the concave case, $z$ decreases).\\

\noindent For the 2-wave there is only a  contact discontinuity (CD):
\begin{equation}
\begin{cases}
\begin{aligned}
&w=w_-\\
&z=z_- +\sigma, \quad \sigma \in \R
\end{aligned}
\end{cases}
\label{aS2}
\end{equation}

\subsection{The Riemann problem}

\begin{figure}[h]
\includegraphics [scale=0.6]{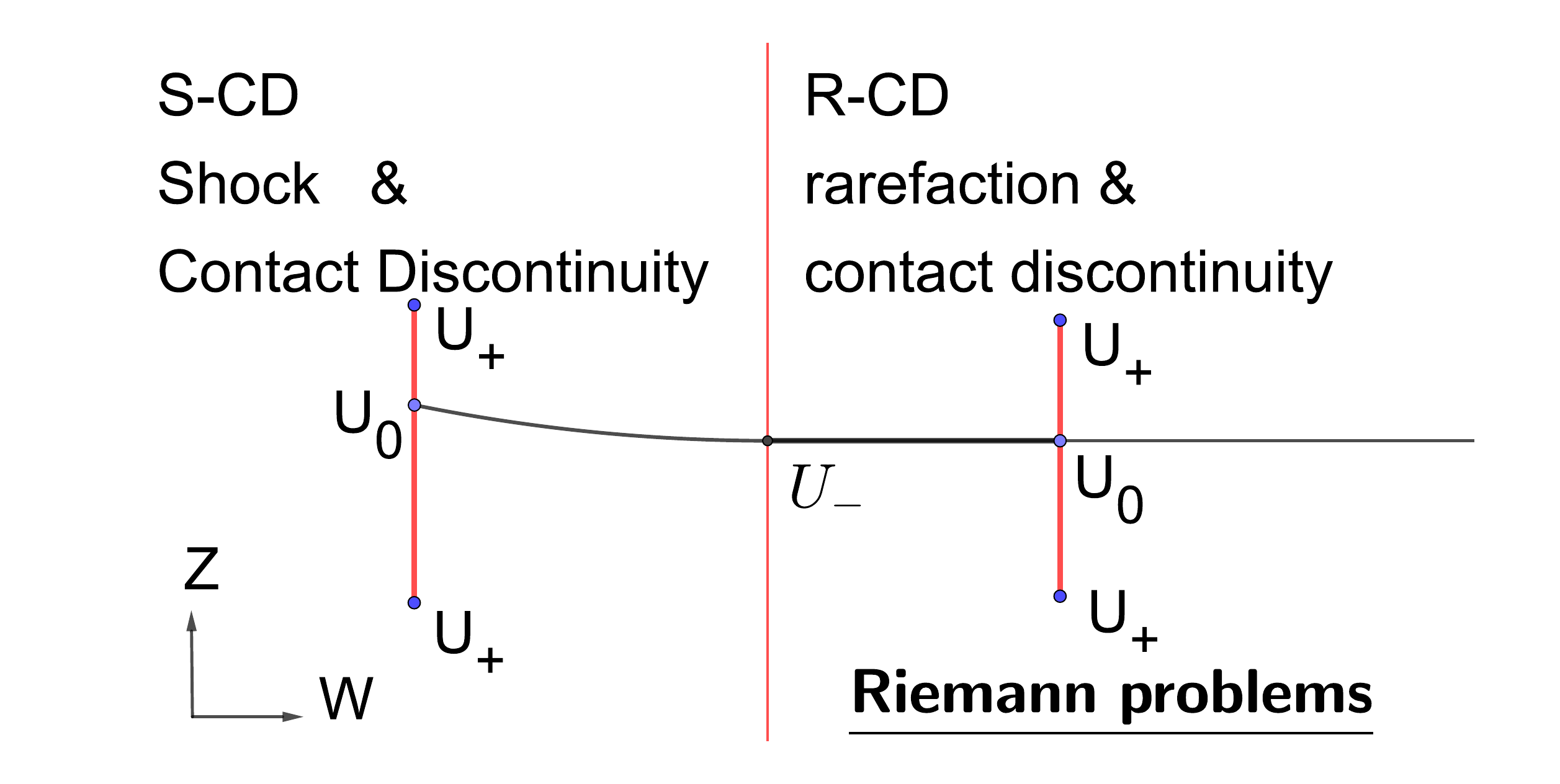}
\caption{Solutions of the Riemann problems for a left state $U_-$ fixed and all the possible configurations for $U_+$. The solutions are represented in the plane $(w,z)$ of Riemann invariants}
\label{fig:Riemann-pb}
\end{figure}

The solution of the Riemann problem is given in the plane $(w,z)$ 
 in figure \ref{fig:Riemann-pb}. 
The initial data is 
$
  U(0,x)=   U_\pm, \quad \pm x \geq 0.
$
$U_0$ is the intermediate  constant state between $U_-$ and $U_+$ when we solve the Riemann problem (do not confuse $U_0$ with the initial data).

\subsection{Nonlinear interactions}
\label{section4.2}

Next we consider the different interactions that we can have. We will note $CD$ for the 2-contact discontinuity wave, $S1$ for the 1-shock wave and $R1$ for the 1 rarefaction wave. We recall that the only possible interactions are:
$$CD-R1,CD-S_1, R_1-S_1, S_1-R_1, S_1-S_1,$$
where $L-R$ means the interaction between a left wave an a right wave. 
The left wave is a 2-wave, a contact discontinuity (CD), or  a 1-wave, a rarefaction (R) or a shock (S) wave. The right wave is always a 1-wave since a 2-wave is a contact discontinuity and goes to the right with the speed $\lambda_2>0$ and cannot interact with a left 1-wave which goes to the left  with the speed $\lambda_1 <0$ or a contact discontinuity which moves with the same speed $\lambda_2$.

We study all the possible interactions in this section. 
The strength of an i-wave is quantified by the variation of associated Riemann invariant through the wave: $\Delta w$ for a 1-wave and  $\Delta z$ for a 2-wave.

 Let us summarize  important features of such interaction  where the 1-wave can be  only a rarefaction or a shock wave and the 2-wave can be only a contact discontinuity. 
  The following list of claims are verified  just after by the exhaustive study of all possible interactions.

\begin{enumerate}
\item When two waves interact  then there are two resulting waves which are a 1-wave on the left and a 2-wave on the right. 
\item  The strength of a 1-wave does not change after an interaction with a 2-wave.
\item  Assume that two 1-wave interacts with respectively the strength $\sigma_1$ and $\widetilde{\sigma}_1$ then the outgoing 1-wave has the strength $\sigma'_1=\sigma_1+\widetilde{\sigma}_1$.
\item  The variation of $w$ after an interaction behaves like the variation of the solution of a scalar conservation law. It means that $TV w$ and $TV^s w$ is not increasing as for a scalar conservation laws.
\item The $L^\infty$ norm of $z$ can increase only when there is an interaction   $CD-S_1$ or $S_1-S_1$.
\end{enumerate}

In all the following pictures,
\begin{itemize}
\item $U_-=(w_-,z_-)$ is the left state, 
\item $U_+=(w_+,z_+)$ is the right state,
\item $U_0=(w_0,z_0)$ is the intermediary state before the interaction,
\item $U_m=(w_m,z_m)$ is the intermediary state after the interaction.
 \end{itemize}

\begin{figure}[t]
\includegraphics [scale=0.8]{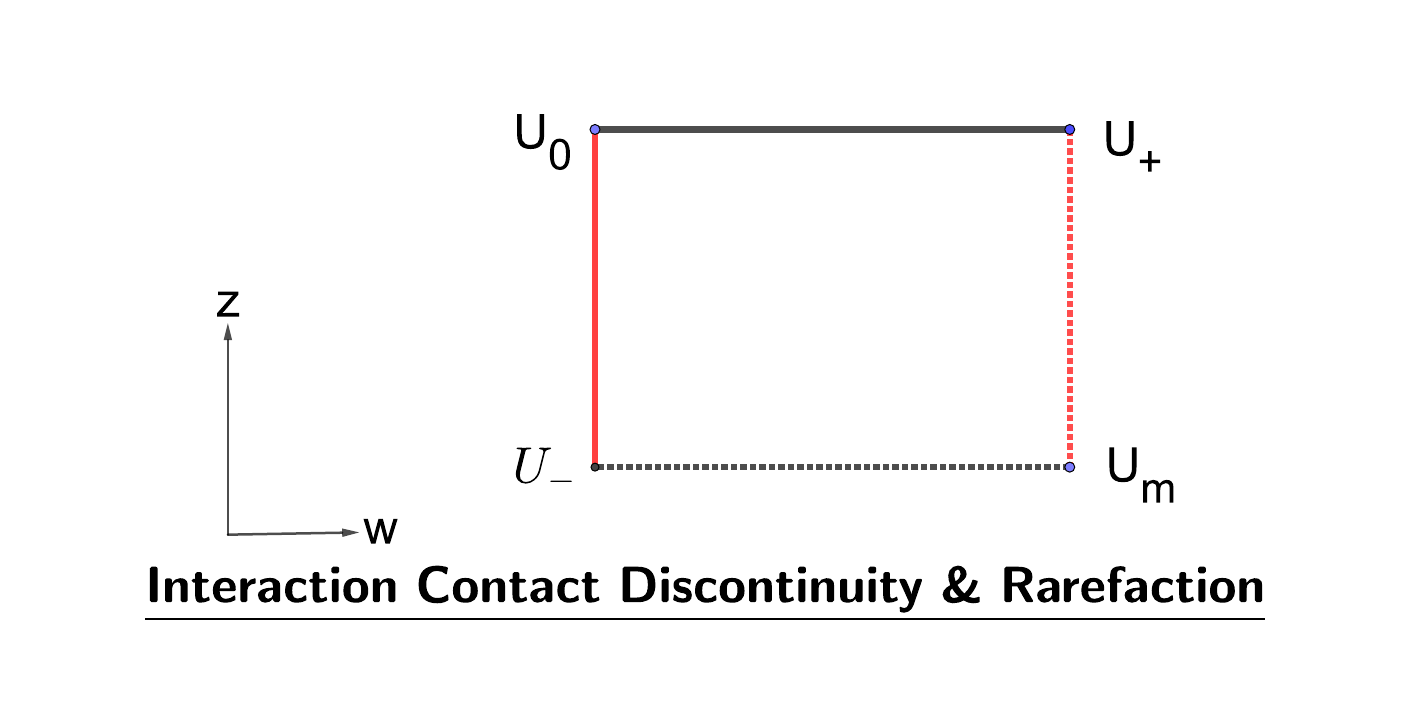}
\caption{Interaction of a contact discontinuity with a rarefaction.  The interacting waves  are represented by FULL lines,  a 2-wave or a 1-wave followed by a 1-wave. The dotted lines represent the resulting waves, a 1-wave (horizontal) followed by a 2-wave (vertical).}
\label{fig:Interaction-DR}
\end{figure}


\begin{figure}[!b]
\includegraphics [scale=0.8]{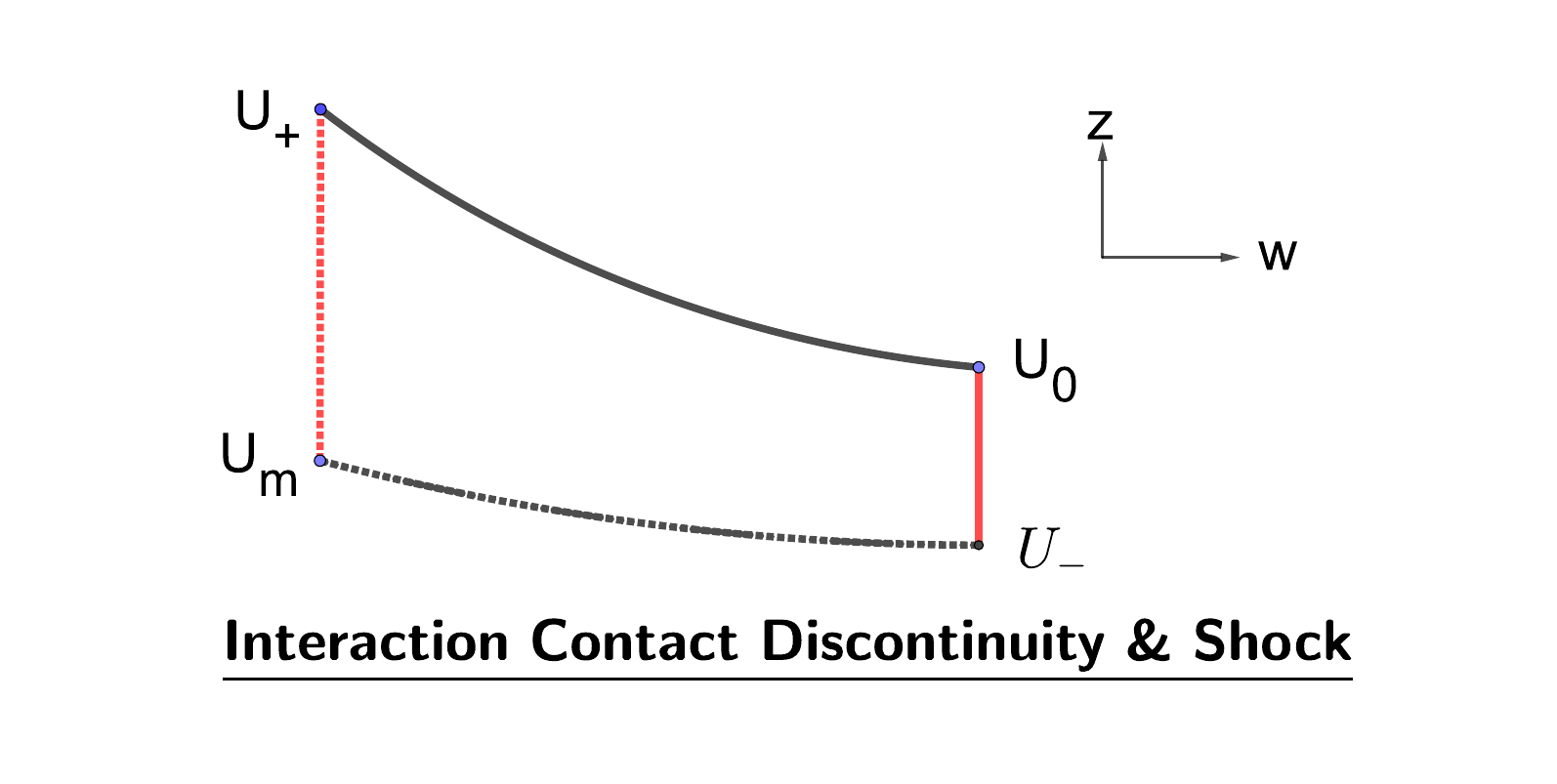}
\caption{Interaction of a contact discontinuity with a shock}
\label{fig:Interaction-DS}
\end{figure}
{ An important point is the control of the $BV$ or $BV^s$ norm of $w$ after an interaction.
 There are two cases.  First case, after the interaction, the solution has only three different values the state $U^-,U_m,U^+$. This is true if  the 1 outgoing wave is a 1-shock.  Second case, there is a  is a 1 rarefaction outgoing  wave, so the solution has a continuum of values (when we will deal with the wave front tracking this continuum of values will be split in a finite number of values depending on the parameter $\nu$ with $\nu$ goes to $+\infty$).  However since the 1 rarefaction wave create a zone of monotonicity in $w$, using the Lemma \ref{lemcru} it does not change all the arguments which are related to estimating the $BV$ or $BV^s$ norm of $w$ (we can apply the Lemma \ref{lemcru}  because the wave front tracking will have a finite number of values).
}
\begin{figure}[!b]
\includegraphics [scale=0.8]{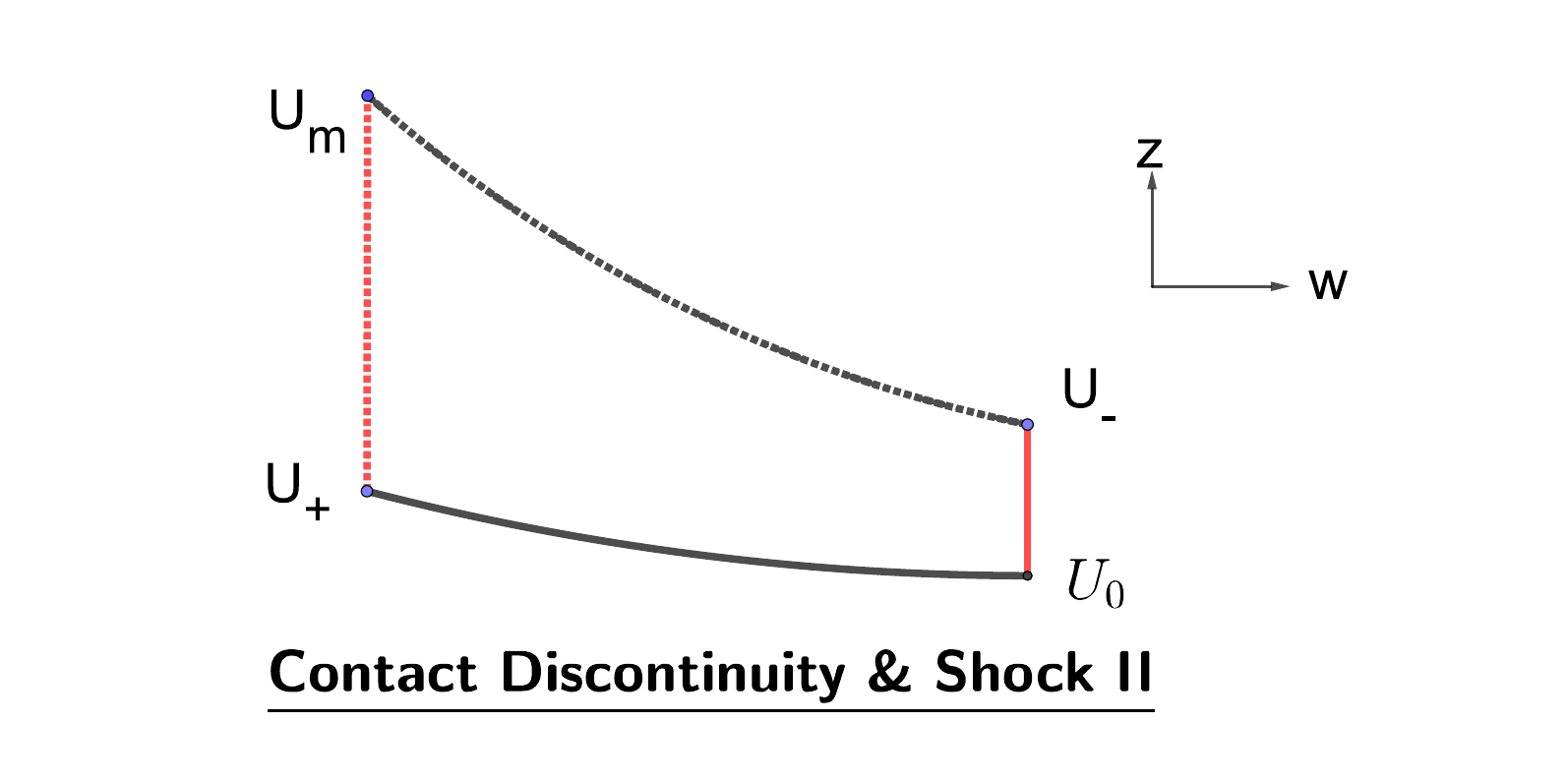}
\caption{Interaction between a contact discontinuity and a shock with an augmentation of $\| z\|_\infty$.}
\label{fig:Interaction-DS2}
\end{figure}

\medskip 
\noindent $\underline{CD-R_1}$. 
 Let us consider the simplest interaction $CD-R_1$, figure \ref{fig:Interaction-DR}, 
which generate waves $R_1-CD$. 
We observe that:
$$w_-=w_0<w_+,\;z_-<z_0=z_+\;\;  \mbox{and}\;\;w_-<w_m=w_+,\;z_-=z_m<z_+.$$
In particular the functions $w(t,\cdot)$ and $z(t,\cdot)$ have the same values before and after the interaction and these values are { in the same order}, it implies then that $BV$ and $BV^s$ norm does not change for this interaction both for $z$ and $w$. Furthermore  we have $w_+-w_0=w_m-w_-$ then the strength of the 1-wave is invariant after this interaction with this 2-wave.

\begin{figure}[t]
\includegraphics [scale=0.6]{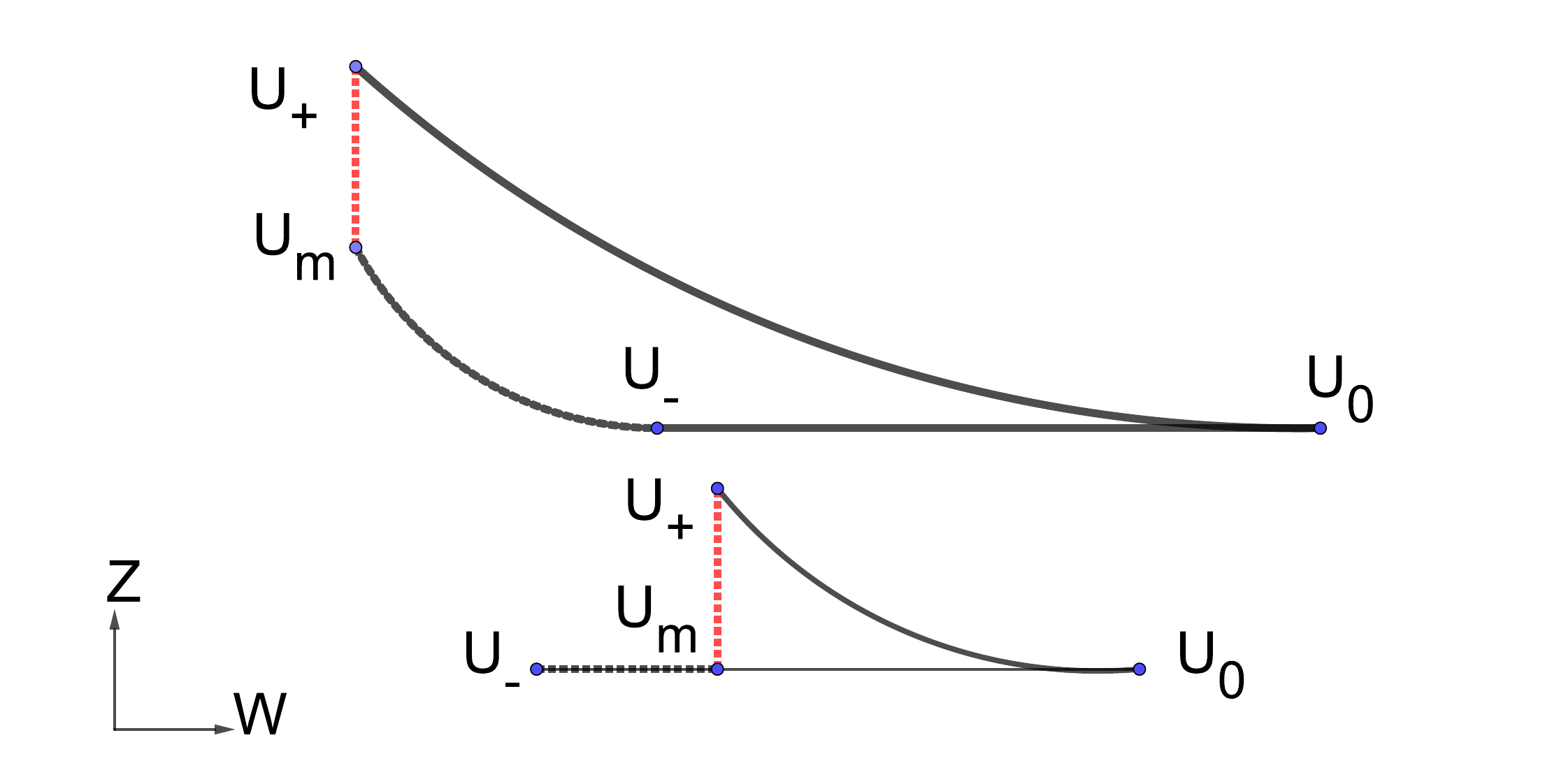}
\caption{Interaction of a rarefaction with a shock wave: two cases depending on the relative strength of the waves}
\label{fig:Interaction-RS}
\end{figure}

\begin{figure}[!h]
\includegraphics [scale=0.6]{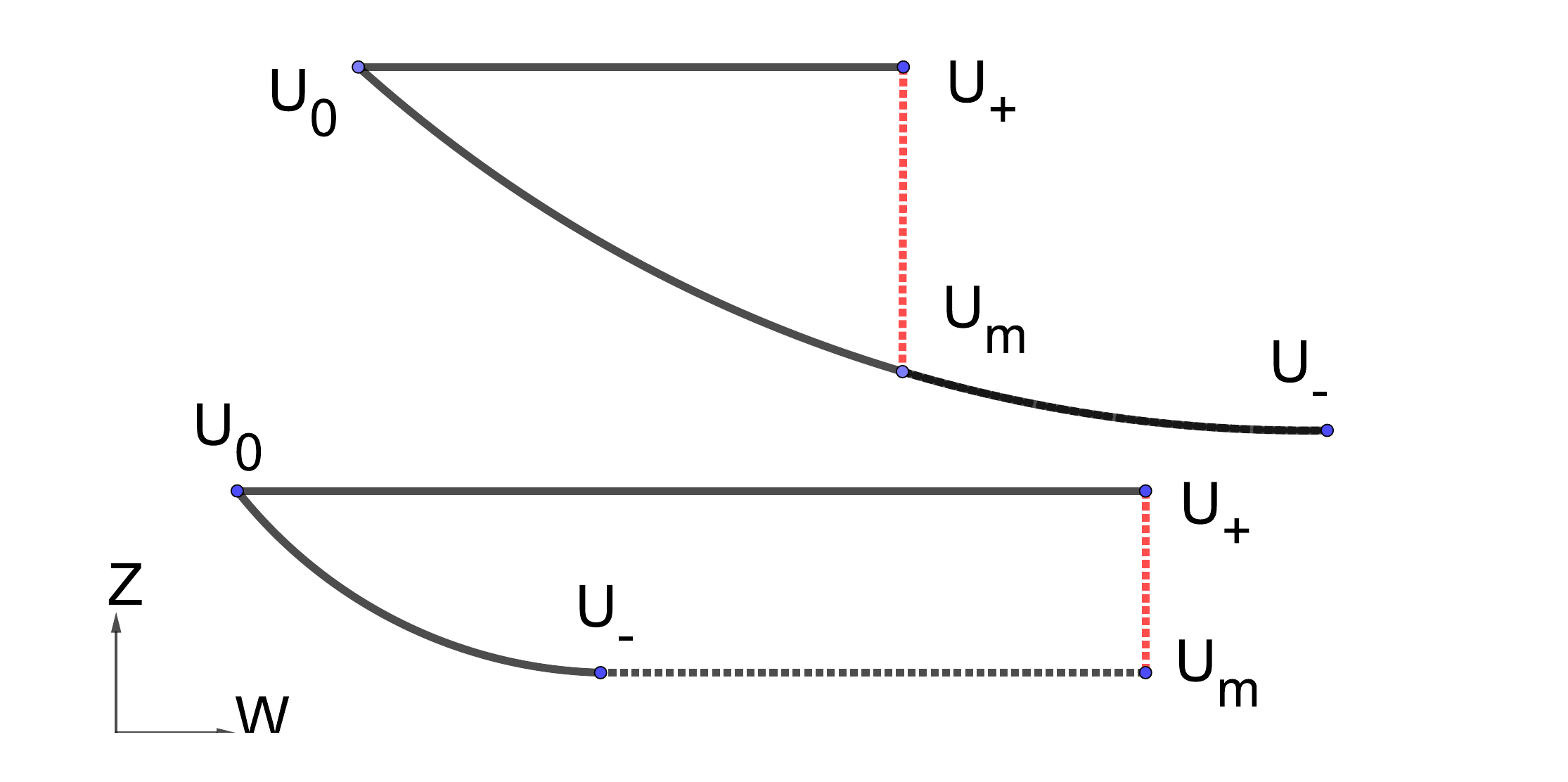}
\caption{Interaction of a shock wave with  a rarefaction: two cases depending on the relative strength of the waves}
\label{fig:Interaction-SR}
\end{figure}

\medskip
\noindent $\underline{CD-S_1}$. Now, consider the interaction $CD-S_1$ , figure \ref{fig:Interaction-DS} and  \ref{fig:Interaction-DS2}, which generates waves $S_1-CD$. We observe that:
$$w_-=w_0>w_+\;\;\mbox{and}\;\;w_0>w_m=w_+ $$
In particular it implies again that the $BV$ and the $BV^s$ norm of $w$ does not change after this type of interaction. Furthermore the strength of the 1 outgoing wave  is the same as the strength of the 1 incoming wave. 
We observe however that the $L^\infty$ norm of $z$ can increase   in the figure \ref{fig:Interaction-DS2},  $|z_m|$ is larger than $|z_-|$, $|z_0|$ and $|z_+|$. Similarly the $BV$ norm can increase for $z$. Thus, there is no maximum principle for $z$.  However, the increase of $\|z\|_\infty$  is controlled as in the last case (\ref{intercru}), the shock-shock interaction.
\medskip


\noindent $\underline{R_1-S_1}$. The interaction of 1-waves $R_1-S_1$ generates $R_1-CD$ or $S_1-CD$,  figure \ref{fig:Interaction-RS},  and we have in each case
$w_m=w_+$. It implies in particular that $w(t,\cdot)$ has the same values after the interaction excepted the value $w_0$ and some values of the incoming 1 rarefaction. \tcb{ Furthermore the values are in the same order, we deduce then that   the $BV$ and the $BV^s$  norm decreases since we restrict in some sense the number of possible subdivision.} 
We can observe that for this interaction the $L^\infty$ norm of $z$ does not increase.
%

\medskip

\noindent $\underline{S_1-R_1}$. The interaction $S_1-R_1$ is similar to the previous case, figure \ref{fig:Interaction-SR}, decay of $TV w$ and invariance of $TV z$.   
\medskip


\begin{figure}[h]
\includegraphics [scale=0.7]{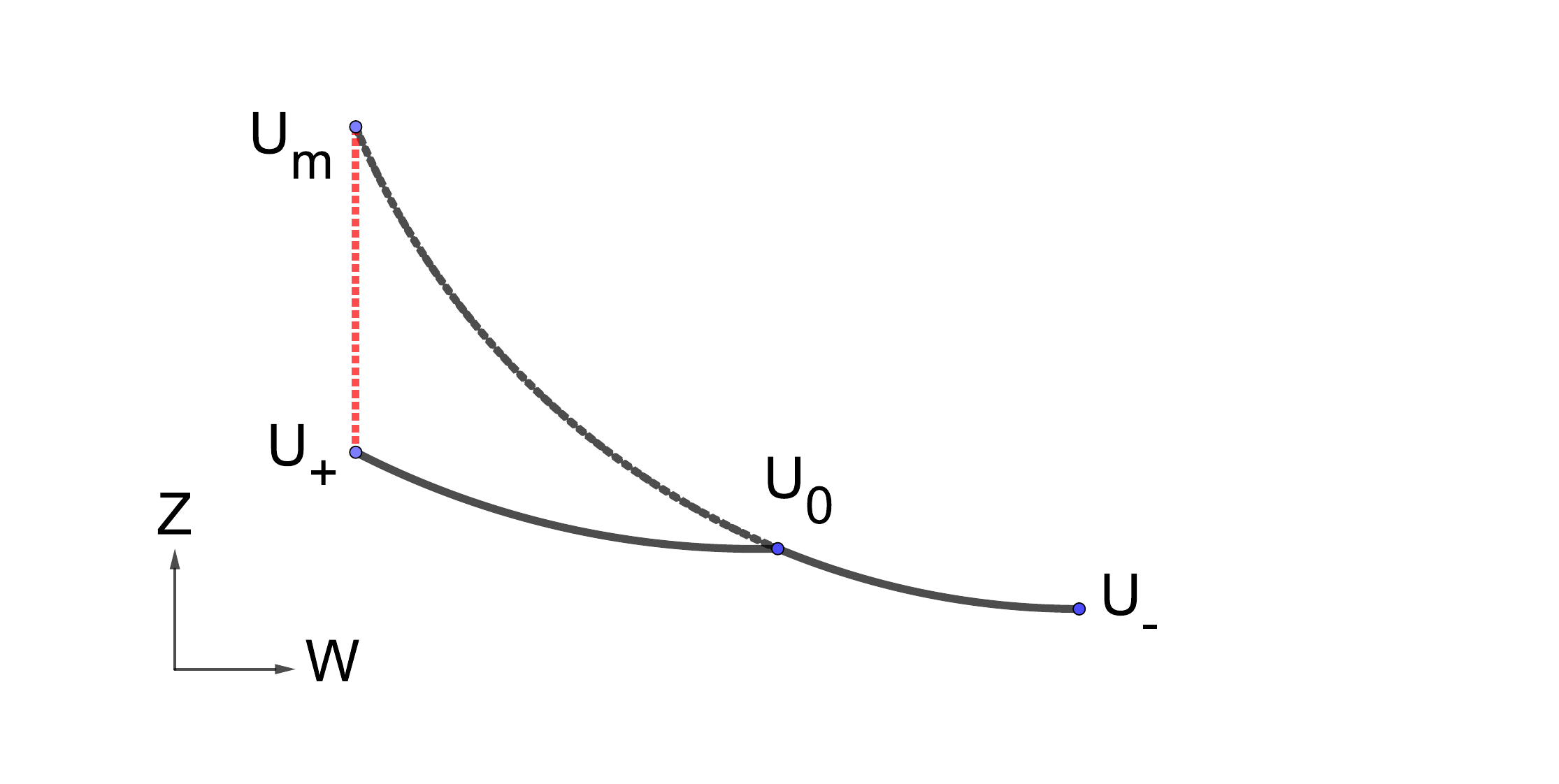}
\caption{Interaction of  two shock waves}
\label{fig:Interaction-SS}
\end{figure}

\noindent $\underline{S_1-S_1}$. We finish with the interaction $S_1-S_1$ which generate the waves $S_1-CD$, figure \ref{fig:Interaction-SS}. This is the most interesting case with a non-scalar type interaction.

Here, $w$ continues to behave  like a solution for a convex scalar conservation law since $w(t,\cdot)$ has the same values after the interaction except $w_0$, { the order of local minimal and local maximal  values of $w$ does not change then the $BV$ and $BV^s$ norms do not  increase for $w$.}
$ z $ is not monotonous after the shock interactions in particular the $L^\infty$ and $BV$ norm can increase. From (\ref{aS1}) one obtain (see also the Glimm estimate in \cite{Da} for a $2 \times 2$ system):
\begin{equation}
\begin{aligned}
   & |z_m-z_+|+|z_m-z_{-}|
   \\
=  & 
|z_+-z_0| + |z_0-z_-|+  2|z_m-z_+|\\
=  & 
|z_+-z_0| + |z_0-z_-|+  2(|z_m-z_-|-|z_0-z_+|-|z_0-z_-|)\\
=  & 
|z_+-z_0| + |z_0-z_-|+  2(O(|w_+-w_-|^3)-O(|w_0-w_+|^3)-O(|w_0-w_-|^ 3))\\
\end{aligned}
\label{intercru}
\end{equation}

In conclusion in all cases,  $w$ always behaves like a solution of a convex scalar law  (in particular the $BV$ and the $BV^s$ norm decreases) and $z$ can only increase (or decrease for concave shock curves) when two shocks interact. In particular the $L^\infty$ and $TV$ norm of $z$ can increase.

%
%

\subsection{A simplified WFT for  such $2\times 2$ systems}
\label{secWFT}
 
%
%
%
%
%


Due to the  special structure of the  $2\times 2$ systems studied,  the  wave front tracking can be simplified in the plane $(w,z)$.
We define now the wave front tracking that we will use in the sequel. First we shall work with initial data $(w_{0,\nu},z_{0,\nu})$ which are piecewise constant approximation of $(w_0,z_0)$ such that:
\begin{equation}
\begin{aligned}
&\|(w_0,z_0)-(w_{0,\nu},z_{0,\nu})\|_{L^1}\leq \nu^{-1},\\
&\mbox{Osc}(w_{0,\nu})\leq \mbox{Osc}(w_0).
\end{aligned}
\label{initial}
\end{equation}
with $\nu>0$ and $\nu$ goes to $+\infty$. Here $\mbox{Osc}(w_0)=\sup_{x,y\in\R}|w_0(x)-w_0(y)|$ denote the oscillation of $w_0$. Furthermore we assume that:
\begin{equation}
w_{0,\nu}\in\nu^{-1}\mathbb{Z}\;\;\mbox{and}\;\;z_{0,\nu}\in\R.
\label{condiini}
\end{equation}
We define now $N_\nu$ as the number of {discontinuities} that the initial data $(w_{0,\nu},z_{0,\nu})$ has. We start the wave front tracking by solving the $N_\nu$ first problem of Riemann. 
We would like to explain how we describe the solution of a Riemann problem between $(w^-,z^-)$ and $(w^+,z^+)$ in our algorithm of wave front tracking (for the beginning we assume that $(w^-,z^-)$ and $(w^+,z^+)$ are some values of $(w_{0,\nu},z_{0,\nu})$). The Riemann problem is the combination of a 1-shock or a 1 rarefaction with a 2-contact discontinuity.
If we get a 1-shock and a contact discontinuity the solution is the exact  solution of the Riemann problem  
with respectively the speed $\lambda_1((z^-,w^-),(z^m,w^m))$, $\lambda_2((z^+,w^+))$ where $\lambda_1((z^-,w^-),(z^m,w^m))$ corresponds to the speed define by the Rankine Hugoniot relation. $(z^m,w^m)$ corresponds here to the intermediary state and we know  that $w^m=w^+$ because the value of $w$ is constant along the 2 discontinuity of contact.
\paragraph{Approximate Riemann solver}
To stay with piecewise constant solutions, an approximate solver is needed only for rarefaction wave.  
If the solution of the Riemann problem is a combination of a 1 rarefaction wave and a 2-contact discontinuity, we have to define the solution corresponding to the rarefaction. We note again $(z^m,w^m)$ the intermediary state of the exact solution. We observe using (\ref{aR1}) and (\ref{aS2}) that $w_m=w_+$, it implies in particular that since $w_+$ is in $\nu^{-1}\mathbb{Z}$ that $w_m$ is again in $\nu^{-1}\mathbb{Z}$. We have in particular:
$$w_m=w_{-}+ k_+ \, \nu^{-1},$$
with $k_+\in\mathbb{N}$ since { $w$ increases through a rarefaction wave}. We define now the intermediary state $w_{k}=w_{-}+ k\, \nu^{-1}$ with $0\leq k\leq k_+$. The solution of the rarefaction for our wave front tracking at time $t>0$ with the initial discontinuity at time $t=0$ in $y$ is:
\begin{equation}
\begin{aligned}
w(t,x)&=w_{-}\;\;\;\mbox{if}\;\;x<x_{1}(t)\\
&=w_{k}\;\;\;\mbox{if}\;\; x_{k}(t)<x<x_{k+1}(t)\\
&=w_{m}\;\;\;\mbox{if}\;\;x>x_{k^+}(t)\\
z(t,x)&=z^-,
\end{aligned}
\label{rar1}
\end{equation}
with $x_{i}(t)=y+t\lambda_1(w_i).$
The approximate solution $u$ can be prolonged until a time $t_1$  when the first interaction between two or more waves front takes place. 
\begin{remarka}
It is important to observe that for $t\in]0,t_1[$, the solution $w(t,x)$ take his values in $\nu^{-1}\mathbb{Z}$.
\end{remarka}
Since $u(t_1,\cdot)$ is still a piecewise constant function, the corresponding Riemann problem can again be approximately solved. The solution $(z,w)$ is then continued up to a time $t_2$ when the second set of wave interactions takes place, and so on.
\begin{remarka}
We assume in the sequel that we have only one interaction at the same time and that the interaction concerns only two fronts and no more (indeed as in \cite{Br}, if three fronts or more meets at the same point, we can avoiding this situation by changing the speed of one of the incoming fronts. Of course this change of speed can be chosen arbitrary small.)
\end{remarka}
We are going to prove in the sequel that there is a finite number of interaction so that we can define the wave front tracking on the time interval $(0,+\infty)$. Assume  for the moment that we can only define the wave-front tracking on a time interval $(0,T^*)$ with $T^*<+\infty$ such that there is an infinite number of wave interactions.   {
 We will prove in fact that  necessary $T^*=+\infty$.}
\begin{remarka}
It is important to verify that $(w,z)(t,\cdot)$ remains bounded in a ball $B(0,r)$ with $r>0$ sufficiently small for all time $t\in(0,T^*)$. Indeed using the Lax Theorem, we can solve the Riemann problem only if the oscillation between two states $(w_{-},z_{-})$ and $(w_+,z_+)$ is sufficiently small. 
\label{Linf} 
\end{remarka}
We define now $N_1(t)$ as:
$$N_1(t)=\mbox{number of 1-wave at the time $t$}.$$
We observe easily that for $t\in]0,t_1[$, we have from (\ref{aR1}):
$$N_1(t)\leq \nu \, N_\nu  \,\mbox{Osc}(w_0)+N_\nu.$$
The first term on the right hand side corresponds to the maximal number of rarefaction and $N_\nu$ to the maximal number of 1-shock. Let us estimate $N_1(t)$ after an interaction at the time $t_k$. We start by recalling that if we have a 1 rarefaction wave by definition of the Riemann problem (see (\ref{rar1})), his strength is necessary of size $\frac{1}{\nu}$.
Assume now that we have an interaction between a 1 rarefaction and a 1-shock then the strength of the 1 rarefaction is $\frac{1}{\nu}$ and the strength of the 1-shock is $\frac{k}{\nu}$ with $k\in-\mathbb{N}^*$ then we have seen in the section \ref{section4.2}  that the 1-wave has the strength $\sigma'_1=\frac{1}{\nu}+\frac{k}{\nu}\leq 0$. It implies in particular that the outgoing 1-wave disappears or is a shock. In particular after such interaction the number of 1-waves decrease of 1 or 2 units. Similarly
if we have an interaction between two shock, we know that the strength of the outgoing 1-wave is the sum of the two strength of the incoming waves, then this strength is negative and the outgoing 1-wave is a shock. It implies in particular that after such interaction the number of 1-waves decrease of 1 unit. If we have an interaction between a 2-wave and a 1-wave since this interaction is transparent the outgoing wave is a unique rarefaction if the incoming one is also a rarefaction and otherwise a 1-shock. It implies in particular that $N_1(t)$ does not change after such interactions.
We have then proved that $N_1(t)$ is a decreasing function of the time and then that:
$$N_1(t)\leq \nu \, N_\nu  \, \mbox{Osc}(w_0)+N_\nu,$$
for all $t\in(0,T^*)$. We define now $N'_1$ as:
$$N'_1=\mbox{Number of interactions between 1-waves on $(0,T^*)$}$$
We have seen that when we have an interaction between 1-waves the number of 1-waves decrease at least of one unit it implies then that $N'_1$ is inferior to the maximal number of 1-waves:
\begin{equation}
N'_1\leq  \nu \,  N_\nu  \, \mbox{Osc}(w_0)+N_\nu.
\label{inter1}
\end{equation}
We define now $N'_2(t)$ as:
$$N'_2=\mbox{Number of interactions between 1-waves and 2-waves on $(0,T^*)$}$$
For a one wave we can define a 1 polygonal line which  is an extension of the one wave. Indeed the one wave are created at the time $t=0^+$ and after each interaction it can be prolongated by a unique 1-wave (or even the one wave can disappear, in this case the polygonal line is stopped). We can then define  a polygonal 1-wave line. We note that two different 1 polygonal lines are similar after an interaction time corresponding to an interaction between 1-waves if they meet us.
Their number is finite and bounded by $ \nu\,  N_\nu \, \mbox{Osc}(w_0)+N_\nu$. Similarly we can define some polygonal 2-wave.  At the difference with the 1 polygonal line, we can create a 2 polygonal line after an interaction between 1-waves. Their number is finite and bounded by $N_\nu+ (N_\nu  \, \mbox{Osc}(w_0)+N_\nu)$, $N_\nu$ corresponds to the number of polygonal line issue of the time $t=0$ and $(N_\nu  \, \mbox{Osc}(w_0)+N_\nu)$ is the maximal number of interaction between 1-waves and then the maximum number of 2 polygonal line that we can create.
Since $\lambda_1(w,z)<0<\lambda_2(w,z)$ for any $(w,z)\in B(0,r)$ we deduce by transversality  that a polygonal 2-wave can interact with a polygonal 1-wave only one time. It implies in particular that the number of interaction on $(0,T^*)$ between polygonal 1-wave and polygonal 2-wave is at more $(2N_\nu+ N_\nu  \, \mbox{Osc}(w_0))( \nu \,  N_\nu \, \mbox{Osc}(w_0)+N_\nu)$. It implies in particular that:
\begin{equation}
N'_2 \leq (2N_\nu+ N_\nu  \, \mbox{Osc}(w_0))( \nu \,  N_\nu \, \mbox{Osc}(w_0)+N_\nu).
\label{inter2}
\end{equation}
From (\ref{inter1}) and (\ref{inter2}), we deduce that the number of interaction on $(0,T^*)$ is finite and then $T^*=+\infty$.
\begin{remarka}{
The only point to verify is to ensure that all along the algorithm of wave front tracking, $(w(t,\cdot),z(t,\cdot))$ must remains in a set $[-r',r']^2$ with $r'>0$ sufficiently small such that we are able to solve any Riemann problem. It will be verified in the sequel when we will prove the Theorem \ref{theo2}. We mention however that the $L^\infty$ norm of $w(t,\cdot)$ is  not increasing. After each interaction, we have proved that the $L^\infty$ norm of $w$ does not increase. It is not the case for $z$ since after the interaction  betwenn two 1-shocks  or the interaction between    a contact discontinuity  and a 1-shock, the $L^\infty$ norm of $z$ can increase.}
\label{remimp}
\end{remarka}




\section{Existence for $s \geq 1/3$  in $BV^s \times L^\infty$}
\label{sec4}
In this section we are going to prove the Theorem \ref{theo2} with the initial data $(w_0,z_0)$ belonging to $BV^s\times L^\infty$ with $s\geq 1/3$.
We consider again the solution of the wave front tracking $(w_\nu,z_\nu)$ defined in the section \ref{secWFT} on a time interval $(0,T_\nu ^*)$ with $T^*_\nu>0$. In addition we construct $(w_{0,\nu},z_{0,\nu})_{\nu>0}$ verifying (\ref{initial}) and such that for any $\nu>0$ we have:
\begin{equation}
\|(w_{0,\nu},z_{0,\nu})\|_{BV^s \times L^\infty}\leq \|(w_{0},z_{0})\|_{BV^s \times L^\infty}.
\label{aini2}
\end{equation}
We are now going to obtain uniform estimate in $\nu$ in $BV^s$ for the solution $w_\nu$. More precisely we wish to prove that
 for any $t\in (0,T_\nu^*)$:
$$\|w_\nu(t,\cdot)\|_{BV^s(\R)}\leq \|w_0\|_{BV^s}.$$

\subsubsection*{Control of $\|w_\nu(t,\cdot)\|_{BV^s}$}
Assume that $t_1\in(0,T_\nu^*)$ is the first time where we have a wave interaction in our wave front tracking. Let us prove now that for any $t\in(0,t_1)$ we have:
\begin{equation}
\|w_\nu(t,\cdot)\|_{BV^s(\R)}=\|w_{0,\nu}\|_{BV^s(\R)}.
\label{TV0}
\end{equation}
If we come back to the estimate (\ref{TV0}), it suffices to observe that the solution $w_\nu(t,\cdot)$ for $t\in(0,t_1)$ is the combination of the solutions of different Riemann problems which deal with all the initial discontinuities. We obtain then a combination of a 1-wave and a 2 CD wave. If we have a 1-shock and a 2 CD wave the values of $w_\nu(t,\cdot)$ does not change compared with $w_\nu(0,\cdot)$ and conserves the same order, then the $BV^s$ norm remains the same. It the 1-wave is a rarefaction, we get different 1 rarefaction fronts and $w_\nu(t,\cdot)$ takes new values. For example if we have a discontinuity in $x_\alpha$ at the time $t=0$ with the values $w_{0,\nu}(x_\alpha^-)$ and $w_{0,\nu}(x_\alpha^+)$ then 
the 1 rarefactions fronts produce the following new values at time $t\in(0,t_1)$ which are $w_{0,\nu}(x_\alpha^-)+\frac{k}{\nu}$ with $k\in\{0,\cdots,k^+\}$ and with $w_{0,\nu}(x_\alpha^+)=w_{0,\nu}(x_\alpha-)+\frac{k^+}{\nu}$.
However even if we have new values for $w_\nu(t,\cdot)$ we have a zone of monotonicity for the 1 rarefaction fronts and using again the Lemma \ref{lemcru}, we conclude again that the $BV^s$ norm does not change. It proves the estimate  (\ref{TV0}).
\\
Next we would like to understand how the $BV^s$ norm vary after each interactions.  Assume that we have an interaction at a time $t_k$ with $U^-,U_0,\, U^+$ the incoming states and  $U^-,U_m,\, U^+$ 
 the outgoing states (here for simplicity of notation we have skip the index $\nu$),  we can observe that the number of different values in $w_\nu$ decrease or remains constant after the interaction. Indeed even when the outcoming 1-wave is one rarefaction, we recall that there is  no more that one rarefaction front (this is due to the statements 2 and 3 p 13), furthermore we have:
 $$U^-=(w^-,z^-),U_m=(w^+,z^m)\;\;\mbox{and}\;\;U^+=(w^+,z^+).$$
If we note $t_k$ the time of interaction, it implies in particular that the value $w_0$ disappears at the time $t_k^+$ and that the values of $w_\nu(t_k^{+},\cdot)$ have the same order as the values of $w_\nu(t_k^{-},\cdot)$. We deduce then that for any interaction, we have:
\begin{equation}
\|w_\nu(t_k^{+},\cdot)\|_{TV^s}\leq \|w_\nu(t_k^{-},\cdot)\|_{TV^s}.
\label{ainSS1}
\end{equation}

 we have then two possibility the outcoming waves are a 1-shock wave and a 2 CD wave or a set of 1 rarefaction fronts and a 2 CD wave.
 In the first case, we observe that $w_\nu(t_k^+,\cdot)$ has the same values as $w_\nu(t_k^-,\cdot)$
 excepted the value $w_0$, furthermore these values have the same order. It implies in particular that:
\begin{equation}
\|w_\nu(t_k^{+},\cdot)\|_{TV^s}\leq \|w_\nu(t_k^{-},\cdot)\|_{TV^s}.
\label{ainSS1}
\end{equation}
We deduce from (\ref{TV0}), (\ref{ainSS1}) and (\ref{aini2}) that the norm $\|w_\nu(t,\cdot)\|_{BV^s(\R)}$ is decreasing all along the time and in particular we have for any $t\in(0,T_\nu^*)$:
\begin{equation}
\|w_\nu(t,\cdot)\|_{BV^s(\R)}\leq\|w_{0,\nu}\|_{BV^s(\R)}\leq \|w_{0}\|_{BV^s(\R)}.
\label{TV0a}
\end{equation}
We are now going to  bound uniformly $z_\nu(t,\cdot)$ in $\nu$ in $L^\infty$ norm for any $t\in(0,T^*_\nu)$.
\subsubsection*{Control of $\|z_\nu(t,\cdot)\|_{L^\infty}$}
We recall in particular that this is important to control the $L^\infty$ norm of $z_\nu$ in order to prove that the wave front is well defined, i.e. $U^\nu$ stays in $\Omega$ (we can solve Riemann problem only if the oscillation of $z_\nu$ and $w_\nu$ are sufficiently small, we refer to the Remark \ref{remimp}). To do this, we define $\gamma^\nu_2(t,x_0)$ the  forward generalized 2-characteristic  (see  \cite{Da}) which is an absolutely continuous
solution of the differential inclusion:
$$
\begin{aligned}
&\frac{d}{dt}\gamma^\nu_2(t,x_0)\in[\min (\lambda_2((w_\nu,z_\nu)(t,\gamma^\nu_2(t,x_0)^+),\lambda_2((w_\nu,z_\nu)(t,\gamma^\nu_2(t,x_0)^-)\\
&\hspace{5cm}, \max (\lambda_2((w_\nu,z_\nu)(t,\gamma^\nu_2(t,x_0)^+),\lambda_2((w,z)(t,
\gamma_2(t,x_0)^-)],
\end{aligned}
$$
and such that $\gamma_2^\nu(0,x_0)=x_0$.
In the sequel, in order to simplify the notation, we just will denote by $\gamma^\nu_2(t)$ a  forward generalized 2-characteristic. Now we are interesting in estimating the $L^\infty$ norm of $z_\nu$ along a  forward generalized 2-characteristic such that $\gamma_2^\nu(0)\ne x_\alpha$ with $x_\alpha$ the points where $(w_{0,\nu},z_{0,\nu})$ is discontinuous. In order to follow the evolution of the $L^\infty$ norm of $z_\nu$ along a  forward generalized 2-characteristic, it is important to understand when the $L^\infty$ norm of $z_\nu$ can vary. It is the case only when  the  forward generalized 2-characteristic meets a 1-wave, a 2-wave or an interaction point. We recall that since $\lambda_2>0$ and $\lambda_1<0$ a  forward generalized 2-characteristic and a 1-wave are necessary transversal. 
First let us assume that $\gamma_2$ meets a 1 rarefaction, from (\ref{aR1}) we know that the value of $z_\nu$ is the same along a 1 rarefaction. In particular it implies that the value of $z_\nu$ on $\gamma_2$ does not change when the  forward generalized 2-characteristic meets a 1 rarefaction wave. Now we assume that the  forward generalized 2-characteristic $\gamma_2^\nu$ meets a 1-shock. Again using (\ref{aS1}) and the figure 1, we observe two things. First we remark that $z_\nu^+>z_\nu^-$, it induces in particular that $z_\nu$ is increasing in this case, second we know that:
$$z_\nu^+-z_\nu^-=0((w_\nu^+-w_\nu^-)^3).$$
Now it is important to mention that a  generalized forward 2-characteristic can not meet a 2 CD wave front. Indeed since the second wave is degenerate, we know that $\lambda_2$ does not depends on $z$. In particular using  (\ref{aS2}) we deduce that $\lambda_2((w_\nu^-,z_\nu^-))=\lambda_2((w_\nu^+,z_\nu^+))$ along a 2 CD, it implies that 2 generalized forward characteristic can not meet a 2 CD wave (indeed they should be parallel to the 2 CD front).\\
Now let us deal with the last case when the  forward generalized 2-characteristics meets an interaction point. The only case is when the interaction point is between a 1 rarefaction front and a 1-shock front, or between two 1-shock fronts. Let us start with the case of two 1-shock fronts, we define then by $\sigma_1$ and $\widetilde{\sigma}_1$ the strength of the 2 incoming wave fronts with $U^-,U_0,U^+$ the incoming states such that for $U^-=(w^-,z^-)$:
\begin{equation}
U_0=(w^-+\sigma_1,z^-+0(\sigma_1^3)),\;U^+=(w^-+\sigma_1+\widetilde{\sigma}_1,z^-+O(\sigma_1^3)+O(\widetilde{\sigma}_1^3))
.
\label{example1}
\end{equation}
 We know that the outcoming intermediary state is:
 $$U_m=(w^-+\sigma'_1,z^m)=(w^-+\sigma_1+\widetilde{\sigma}_1,z^-+O((\sigma_1+\widetilde{\sigma}_1)^3).$$
 Now since $\lambda_2$ does not depend on $z$, we deduce that the  forward 2-characteristic follow the outcoming 2 CD wave after the interaction. 
 \begin{remarka}
 By convention, we assume that the value of $z_\nu$ on the 2 CD wave front corresponds to $z_\nu^+$ the value on the right of the 2 CD wave front.
 \end{remarka}
 It implies then that after the interaction point the value of $z_\nu$ has increased on $\gamma_2^\nu$ and is such that:
 \begin{equation}
 z^+_\nu-z^-_\nu=O(\sigma_1^3)+O(\widetilde{\sigma}_1^3).
 \label{bbS1}
 \end{equation}
 Here $z^+_\nu$ is the value of $z_\nu$ on $\gamma_2^\nu$ just after the interaction point and $z^-_\nu$  the value of $z_\nu$ just before.
 Let us consider the case now of an interaction between a 1 rarefaction front and a 1-shock front (the case of an interaction between a 1-shock front and a 1 rarefaction front is similar), we have then:
 \begin{equation}
 U^-=(w^-,z^-),\;U_0=(w^-+\sigma_1,z^-),\;U^+=(w^-+\sigma_1+\widetilde{\sigma}_1,z^-+O(\widetilde{\sigma}_1^3)).
 \label{a6.50}
 \end{equation}
 We know that the outcoming intermediary state is:
 \begin{equation}
 U_m=(w^-+\sigma'_1,z^m)=(w^-+\sigma_1+\widetilde{\sigma}_1,z^-+O((\sigma_1+\widetilde{\sigma}_1)^3)).
 \label{interm2}
 \end{equation}
Again we deduce that the  forward 2-characteristic follow the outcoming 2 CD wave front after the interaction. And we have in addition: 
\begin{equation}
z^+_\nu-z^-_\nu=0(\widetilde{\sigma}_1^3).
\label{bS1}
\end{equation}
 In conclusion we have seen that $z_\nu$ is increasing when the  forward 2-characteristic  meets an interaction point and that we have (\ref{bbS1}) or (\ref{bS1}).
 \\
 Assume now that $\gamma^\nu_2(0)=x_\alpha$, then $\gamma_2^\nu$ is the 2-wave polygonal front which is issue from $x_\alpha$. If $\gamma_2^\nu$ meets 1 rarefaction then the value of $z_\nu$ does not change along the  forward generalized 2-characteristic (see figure 3). If the interaction is with a 1-shock front (see figure 4), we observe that:
 $$z_\nu^+=(z_\nu^+-z^\nu_0)+z^\nu_0=O(\sigma_1^3)+z^\nu_0.$$
 Here $z^\nu_0$ is the intermediary state before the interaction and $\sigma_1$ is the strength of the 1-shock wave. We recall that the value of $z_\nu$ on $\gamma_2$ before the interaction is by convention $z^0_\nu$ since we consider the value on the right for a 2  CD wave front.
 We observe then that if $t_k$ is the interaction point we have:
 $$z_\nu(t_k^+,\gamma^\nu_2(t_k^+))=z(t_k^-,\gamma^\nu_2(t_k^-))+O(\sigma_1^3).$$
 \begin{remarka}
If $\gamma_2^\nu(0)\ne x_\alpha$ and that the 2 forward $\gamma_2^\nu$ becomes after meeting an interaction point a 2 CD polygonal wave front, we can estimate the evolution of the $L^\infty$ norm of $z_\nu$ as in the case where $\gamma_2^\nu(0)= x_\alpha$.
\end{remarka}
  We can now calculate the value of $z_\nu$ at the point $(T,\gamma^\nu_2(T))$ with $T>0$. We have seen using the fact that $z_\nu$ is increasing along $\gamma^\nu_2(t)$ and that $z_\nu$ increases of $O(\sigma_1^3)$ after each interaction with a 1-shock front or an interaction point where there is a 1-shock front, we obtain  then from (\ref{bbS1}) and (\ref{bS1}):
\begin{equation}
z_\nu(T,\gamma^\nu_2(T))=z_\nu(0,\gamma^\nu_2(0))+\sum_{\alpha\in J}O(\sigma_\alpha^3).
\label{cru}
\end{equation}
 Here $J$ corresponds to the set of 1-shock wave fronts which have meet $\gamma^\nu_2$ on the time interval $[0,T]$ including the interaction points. 
 In particular it exists $C>0$ independent on $\nu$ such that:
\begin{equation}
 \sum_{\alpha\in J}O(\sigma_\alpha^3)\leq C\; TV^s w_\nu(\cdot,\gamma_{2}^\nu(\cdot))([0,T])+ O \left(\nu^{-3} \right),
 \label{estimaop}
 \end{equation}
 with $s=\frac{1}{3}$ (we deal always in the sequel with $s=\frac{1}{3}$). 
  \begin{remarka}
 It is important to point out that if  $\gamma_2^\nu(0)\ne x_\alpha$ then the forward 2-characteristic $\gamma^\nu_2$ can meet only one time an interaction point. Indeed after this interaction the forward 2-characteristic becomes a 2 polygonal line and we have constructed a wave front tracking where the interactions concern only two fronts. It means that the forward 2-characteristic after the meeting with an interaction point can cross after only 1 rarefaction fronts and 1-shock fronts. Similarly if $\gamma_2^\nu(0)= x_\alpha$ then the forward 2-characteristic $\gamma^\nu_2$  which is a 2 polygonal front will meet only 1 rarefaction fronts and 1-shock fronts. 
 \label{remtimp}
  \end{remarka}
In order to prove  (\ref{estimaop}), we only consider the case where $\gamma_2^\nu(0)\ne x_\alpha$ and the case where the forward generalized 2-characteristic $\gamma^\nu_2$ meets one interaction point. The other case are simple to treat. Let us start with the case where $\gamma_2^\nu$ meets an interaction point with two 1-shock wave fronts at the time $t_k$, from (\ref{example1}), (\ref{bbS1}) we deduce that:
$$z_\nu(t_k^+,\gamma^\nu_2(t_k^+))=z(t_k^-,\gamma^\nu_2(t_k^-))+O(\sigma_1^3)+0(\widetilde{\sigma}_1^3)\leq z(t_k^-,\gamma^\nu_2(t_k^-))+O(|\sigma_1+\widetilde{\sigma}_1|^3).$$
 And in particular it says that:
\begin{equation}
z_\nu(t_k^+,\gamma^\nu_2(t_k^+))\leq z(t_k^-,\gamma^\nu_2(t_k^-))+O(|w_\nu(t_k^+,\gamma^\nu_2(t_k^+))-w_\nu(t_k^-,\gamma^\nu_2(t_k^-)|^3).
\label{tresimp}
\end{equation}
 Let us deal now with the more tricky case of the cross of $\gamma_2^\nu$ with  an interaction point comprising 1 rarefaction wave front and 1-shock front at the time $t'_k$, we have obtained from (\ref{a6.50}), (\ref{bS1}) that:
\begin{equation}
z_\nu((t'_k)^+,\gamma^\nu_2((t'_k)^+))\leq z((t'_k)^-,\gamma^\nu_2((t'_k)^-))+O(\widetilde {\sigma}_1^3).
\label{hyperimp}
\end{equation}
 We recall that we have:
\begin{equation}
\widetilde{\sigma}_1=w_\nu((t'_k)^+,\gamma^\nu_2((t'_k)^+))-w_\nu((t'_k)^-,\gamma^\nu_2((t'_k)^-)-\sigma_1\leq 0.
\label{m1}
\end{equation}
 It is important to note that $\sigma_1=\frac{1}{\nu}$, indeed our wave front tracking ensures that all the rarefaction fronts have the strength $\frac{1}{\nu}$. Similarly we know from (\ref{interm2}) that the intermediary outcoming state is $U_m=(w_+,z_m)$ and that the 1 outcoming wave is a 1-shock or is cancelled out. For the moment assume that the 1 outcoming wave front is not cancelled out,  in particular since $U_m=(w_+,z_m)$, it implies that $w_+<w_-$. But we know that $w^\nu(t'_k,\cdot)$ takes only values in $\frac{\mathbb{Z}}{\nu}$, we deduce then that:
\begin{equation}
w_\nu((t'_k)^+,\gamma^\nu_2((t'_k)^+))-w_\nu((t'_k)^-,\gamma^\nu_2((t'_k)^-)=-\frac{k}{\nu},
\label{m2}
\end{equation}
 with $k\in\mathbb{N}^*$. From (\ref{m1}) and (\ref{m2}), we have:
\begin{equation}
|\widetilde{\sigma}|^3=\left(\frac {k+1}{\nu}\right)^3\leq 2^3 \left(\frac{k}{\nu} \right)^3=2^3 |w_\nu((t'_k)^+,\gamma^\nu_2((t'_k)^+))-w_\nu((t'_k)^-,\gamma^\nu_2((t'_k)^-)|^3.
\label{hyperimp1}
\end{equation}
 We deduce then using (\ref{hyperimp}) and (\ref{hyperimp1}) that:
\begin{equation}
\begin{aligned}
&z_\nu((t'_k)^+,\gamma^\nu_2((t'_k)^+))\leq z((t'_k)^-,\gamma^\nu_2((t'_k)^-))\\
&\hspace{3cm}+O( |w_\nu((t'_k)^+,\gamma^\nu_2((t'_k)^+))-w_\nu((t'_k)^-,\gamma^\nu_2((t'_k)^-)|^3).
\end{aligned}
\label{hyperimp3}
\end{equation} 
We finish now with the case where the 1 outcoming wave is cancelled out, it corresponds to the following situation:
 \begin{equation}
 U^-=(w^-,z^-),\;U_0=(w^-+\sigma_1,z^-),\;U^+=\left(w^-,z^-+O\left(\sigma_1^3\right)\right).
 \label{b6.50}
 \end{equation}
In this case since $\sigma_1=\frac{1}{\nu}$ because this is the strength of a 1 rarefaction front, we have:
\begin{equation}
\begin{aligned}
&z_\nu((t'_k)^+,\gamma^\nu_2((t'_k)^+))\leq z((t'_k)^-,\gamma^\nu_2((t'_k)^-))+O\left(\nu^{-3}\right).
\end{aligned}
\label{hyperimp4}
\end{equation} 
We proceed similarly for an interaction between a 1-shock wave and a 1 rarefaction wave. 
From the Remark \ref{remtimp}, we know that the forward generalized 2-characteristic $\gamma_2^\nu$ can meet only one time an interaction point then combining (\ref{tresimp}), (\ref{hyperimp3}) and (\ref{hyperimp4})  allows to prove the estimate (\ref{estimaop}).\\
\\
 We wish now to estimate
$TV^s w_\nu(\cdot,\gamma_{2}^\nu(\cdot))([0,T])$ 
  in terms of the $BV^{\frac{1}{3}}$ norm of $w_\nu(0,\cdot)$. To do this we are going to consider a zone of dependence of the  forward generalized 2-characteristic $\gamma^\nu_2$. We now choose the 1 
 minimal backward generalized characteristic issue from $(\gamma^\nu_2(T),T)$ that we note $\gamma^\nu_1$ (note that $\gamma^\nu_1$ is defined on $[0,T]$, we refer to \cite{Da} Chapter X for the notion of minimal backward generalized characteristics).   We now define  the sequence of following functions $\gamma^\nu_{2,\alpha}$ with $\alpha\in[0,1]$ as follows with $\tau\in[0,T]$ (here $\tau$ does not correspond to the physical time $t$), 
 $$\gamma^\nu_{2,\alpha}(\tau)=\alpha\gamma^\nu_2(\tau)+(1-\alpha)\gamma^\nu_2(0)+C^\nu_\alpha \tau,$$
 with 
 $$ C^\nu_\alpha=\frac{\gamma^\nu_1(\alpha T)-\alpha\gamma^\nu_2(T)-(1-\alpha)\gamma^\nu_2(0)}{T}\geq 0$$
  since $\gamma^\nu_1(\alpha T)\geq \gamma^\nu_1(T)=\gamma^\nu_2(T)\geq \alpha\gamma^\nu_2(T)+(1-\alpha)\gamma^\nu_2(0)$ for any $\alpha\in[0,1]$ (indeed we recall that the backward generalized characteristic goes on the right since $\lambda_1<0$). We observe also that:
 $$\gamma^\nu_{2,\alpha}(T)=\gamma^\nu_1(\alpha T)\;\;\mbox{and}\;\;(\gamma^\nu_{2,\alpha})'(t)\geq 0\;\;\mbox{for}\;t\in[0,T].$$
 The derivative of $\gamma^\nu_{2,\alpha}$ is in fact defined on the point where $\gamma_2^\nu$ is differentiable. 
 \begin{remarka}
 It is important to note that the  forward  generalized 2-characteristic is defined in a unique way. The second point is that the domain  delimited by the curves $\{(\gamma^\nu_2(t),t),t\in[0,T]\}$, $\{(\gamma^\nu_1(t),t),t\in[0,T]\}$ and 
 $\{(y,0),y\in [\gamma_2 ^\nu(0),\gamma^\nu_1(0)]\}$ is the union of all the curves $\{(\gamma^\nu_{2,\alpha}(t),\alpha t),t\in[0,T]\}$ with $\alpha\in [0,1]$. We denote by $\Gamma_2$ this domain.
 There is no 1-wave front which enters in $\Gamma_2$ on the right since we have taken the minimal backward 1-characteristic. 
 \end{remarka}
 \begin{figure}[!h]
\includegraphics [scale=0.7]{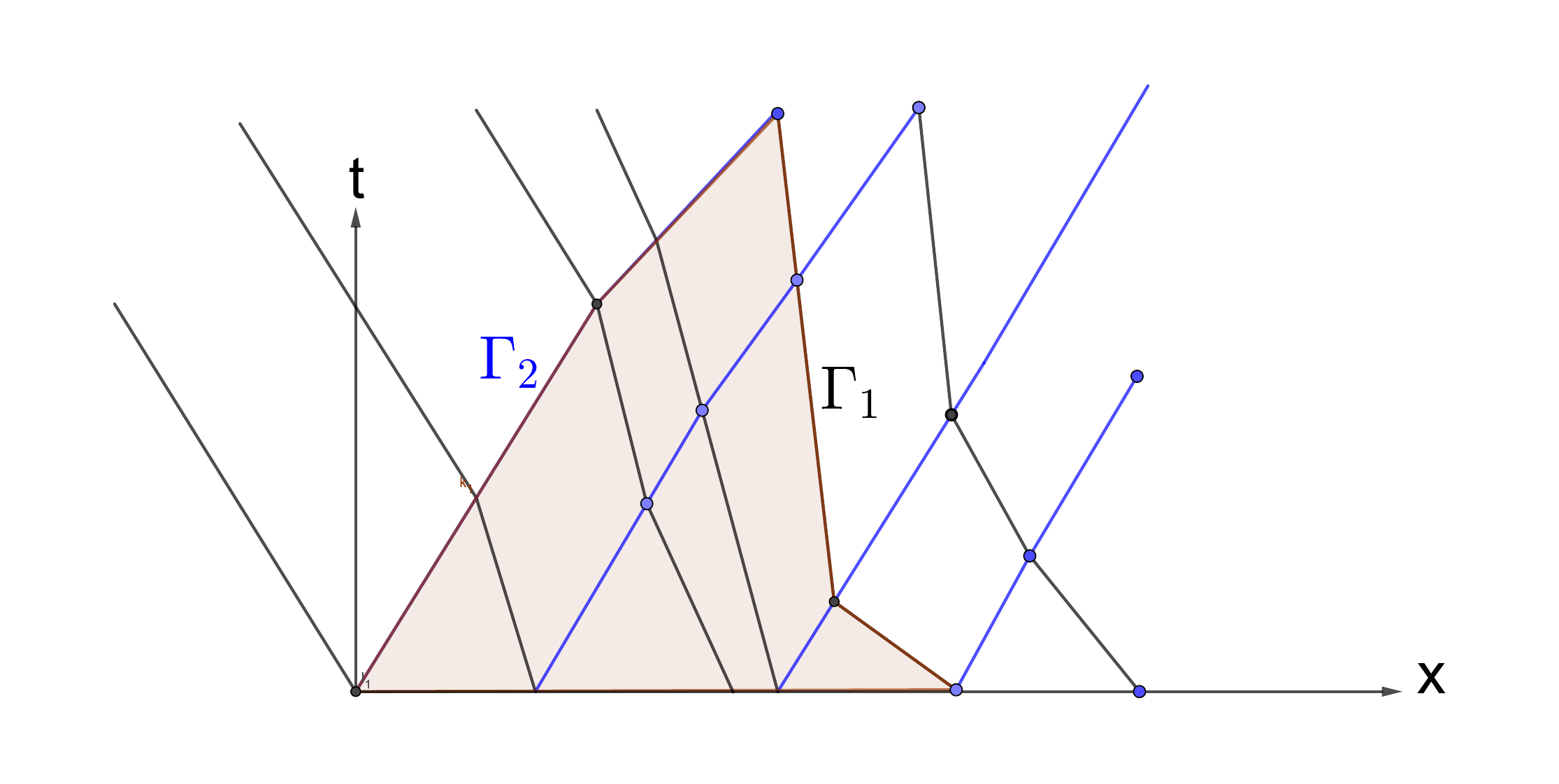}
\caption{The Wave Front Tracking and the dependance zone delimited by the 2-characteristic $\Gamma_2$ on the left, and the 1-characteristic $\Gamma_1$ on the right}
\label{fig:WFT}
\end{figure}
We define now $(x_1,t_1)$ as the point where there is for the first time an interaction between wave fronts inside the domain $\Gamma_2$ with $t_1\in]0,T]$ . We denote now by $t_{1,1}$ the first time where there is an interaction in our wave front tracking, it implies in particular that $0<t_{1,1}\leq t_1$. Furthermore we know that for $\alpha$ such that $\alpha T<t_1$
there is no interaction point on the curve  $\{(\gamma^\nu_{2,\alpha}(t),\alpha t),t\in[0,T]\}$. We denote now by $\alpha_0$ the first $\alpha$ where the curve $\{(\gamma^\nu_{2,\alpha}(t),\alpha t),t\in[0,T]\}$
meets an interaction point $(x_k,t_k)$ inside $\Gamma_2$. In particular it implies that there is no interaction point in the open domain $\Gamma_0$ delimited by $\{(\gamma^\nu_{2,\alpha_0}(t),\alpha_0 t),t\in[0,T]\}$, $\{(\gamma^\nu_1(t),t),t\in[0,T]\}$ and 
 $\{(y,0),y\in [\gamma_2 ^\nu(0),\gamma^\nu_1(0)]\}$. We know that the values of $w$ change on the curve $\{(\gamma^\nu_{2,\alpha}(t),\alpha t),t\in[0,T]\}$ only when this curve meets a 1-wave front.
 Furthermore the only 1-wave front which can cross $\{(\gamma^\nu_{2,\alpha}(t),\alpha t),t\in[0,T]\}$ for $0<\alpha<\alpha_0$ are the 1-wave front which are issue of the set $\{(y,0),y\in [\gamma_2 ^\nu(0),\gamma^\nu_1(0)]\}$ (indeed there is no interaction point in the open domain $\Gamma_0$ and the 1-wave which are outside from $\Gamma_0$ can not enter in $\Gamma_0$), it implies then since the curve $\{(\gamma^\nu_{2,\alpha}(t),\alpha t),t\in[0,T]\}$ are transversal to the 1-wave front that  the values of $w$ on $\{(\gamma^\nu_{2,\alpha}(t),\alpha t),t\in[0,T]\}$ for $\alpha\in]0,\alpha_0[$ are included  in the set of the value of $w_\nu(t,\cdot)$ at the time $t=0^+$ (or in other word $t\in(0,t_{1,1}$) and that they are ranged in the same order. In particular it implies using  (\ref{TV0a}) that for any $\alpha\in(0,\alpha_0)$ we have:
\begin{equation}
 TV^s w_\nu(\alpha \cdot,\gamma^\nu_{2,\alpha}(\cdot))([0,T])\leq
TV^s w_\nu(t ,\cdot)(\R) \leq  TV^s w_{0,\nu}(\R)\leq \|w_0\|^{\frac{1}{s}}_{BV^s},
\label{interac1}
\end{equation}
with $t\in(0,t_{1,1})$.
\\
Next we wish to estimate the $BV^s$ norm of $w^\nu$ along $\{(\gamma_{2,\alpha}(t),\alpha t),t\in[0,T]\}$ for $\alpha=\alpha_0^+$. 
We note that the interaction point $(x_k,t_k)$ on $\{(\gamma_{2,\alpha_0}(t),\alpha_0 t),t\in[0,T]\}$ is by definition inside $\Gamma_2$.
In any case of interaction, if we have incoming states $(U^-,U_0,U^+)$ and outgoing states $(U^-,U_m,U^+)$ then the values of $w^\nu$ around $(x_k,t_k^-)$ are $(w_{-},w_0,w_{+})$ and the the values of $w^\nu$ around $(x_k,t_k^+)$ are $(w_{-},w_{+})$. It means that there is one value $w_\nu$  in less $w_0$ on the curves $\{(\gamma^\nu_{2,\alpha}(t),\alpha t),t\in[0,T]\}$ for $\alpha=\alpha_0^+$ compared with the values of $w_\nu$ on the curves  $\{(\gamma^\nu_{2,\alpha}(t),\alpha t),t\in[0,T]\}$  for $\alpha=\alpha_0^+$.  Furthermore the order of the values of $w_\nu$ on the curves $\{(\gamma^\nu_{2,\alpha}(t),\alpha t),t\in[0,T]\}$ with $\alpha=\alpha_0^+$ and on the curves $\{(\gamma^\nu_{2,\alpha}(t),\alpha t),t\in[0,T]\}$ with $\alpha=\alpha_0^-$ does not change.
It implies that the $BV^s$ norm is decreasing after the interaction along the curves $\{(\gamma^\nu_{2,\alpha}(t),\alpha t),t\in[0,T]\}$ with $\alpha-\alpha_0>0$ sufficiently small. It gives then using (\ref{interac1})
that:
\begin{equation}
 TV^s w_\nu(\alpha_0^+ \cdot,\gamma^\nu_{2,\alpha_0^+}(\cdot))([0,T])\leq
TV^s w_\nu(\alpha_0^- \cdot,\gamma^\nu_{2,\alpha_0^-}(\cdot))([0,T])\leq \|w_0\|^{\frac{1}{s}}_{BV^s},
\label{interac2}
\end{equation}
with $t\in(0,t_{1,1})$.
The previous argument is again true if there is more than 1 interaction on the curve $\{(\alpha_0 t,\gamma^\nu_{2,\alpha_0}(t)),t\in[0,T]\}$ . Now we define $\alpha_1>\alpha_0$ the next $\alpha$ where there is an interaction inside $\Gamma_2$ on the curve $\{(\alpha_1 t,\gamma^\nu_{2,\alpha_1}(t)),t\in[0,T]\}$ and we define by $\Gamma_1$ the open domain delimited by the curves $\{(\alpha_1 t,\gamma^\nu_{2,\alpha_1}(t)),t\in[0,T]\}$, $\{(\alpha_0 t,\gamma^\nu_{2,\alpha_0}(t)),t\in[0,T]\}$ and 
$\{(\gamma^\nu_1(t),t),t\in[0,T]\}$. We observe then that all the 1-wave which meet a curve 
$\{(\alpha t,\gamma^\nu_{2,\alpha}(t)),t\in[0,T]\}$ with $\alpha\in(\alpha_0,\alpha_1)$ are issue of the curve $\{(\alpha_0 t,\gamma^\nu_{2,\alpha_0}(t)),t\in[0,T]\}$. It implies that the values of $w_\nu$ on $\{(\alpha t,\gamma^\nu_{2,\alpha}(t)),t\in[0,T]\}$ are included in set of values of $w_\nu$ on the curve $\{(\alpha' t,\gamma^\nu_{2,\alpha'}(t)),t\in[0,T]\}$ with $\alpha'=\alpha_0^+$, furthermore by transversality the values keep the same order. We deduce then that for any $\alpha\in(\alpha_0,\alpha_1)$ we have using in addition (\ref{interac2}):
\begin{equation}
 TV^s w_\nu(\alpha \cdot,\gamma^\nu_{2,\alpha}(\cdot))([0,T])\leq
 TV^s w_\nu(\alpha_0^+ \cdot,\gamma^\nu_{2,\alpha_0^+}(\cdot))([0,T])\leq \ \|w_0\|^{\frac{1}{s}}_{BV^s}.
\label{interac3}
\end{equation}
Repeating the argument, we deduce finally  that the function:
\begin{equation}
\alpha\rightarrow  \|w(\alpha \cdot,\gamma^\nu_{2,\alpha}(\cdot))\|_{BV^s([0,T])}
\end{equation}
is decreasing in $\alpha$. It implies from (\ref{estimaop}) and (\ref{aini2}) since $\gamma_2^\nu=\gamma_{2,\alpha}^\nu$ with $\alpha=1$ that it exists $C>0$ such that for any $\nu>0$ we have:
\begin{equation}
 \sum_{\alpha\in J}O(\sigma_\alpha^3)\leq C \|w_{0}\|^{\frac{1}{s}}_{BV^s(\R)}+O\left(\nu^{-3}\right),
 \label{crua}
 \end{equation}
 with $s=\frac{1}{3}$. 
From (\ref{cru}) and (\ref{crua}), we deduce that it exists $C>0$ independent on $\nu$ such that for any $T\in (0,T^*_\nu)$ and any forward generalized 2-characteristic $\gamma^\nu_2$:
\begin{equation}
|z_\nu(T,\gamma^\nu_2(T))|\leq |z_\nu(0,\gamma^\nu_2(0))|+C\|w_{0}\|_{BV^s(\R)}^{\frac{1}{s}}+ O\left(\nu^{-3}\right). \label{crub}
\end{equation}
Since the  forward generalized 2-characteristics describe all the space $(0,T_\nu^*)\times\R$, we deduce from (\ref{crub}) and (\ref{aini2}) that for any $t\in(0,T_\nu^*)$ we get for $C>0$ independent on $\nu$:
$$\|z_\nu(t,\cdot)\|_{L^\infty}\leq \|z_0\|_{L^\infty}+C \|w_0\|_{BV^{\frac{1}{3}}(\R)}^3+o(1). $$
We deduce now that the $L^\infty$ norm of $z_\nu$ is uniformly bounded in $\nu$ all along the time interval $(0,T_\nu^*)$ and remains small for large $\nu$, then using the Remark \ref{remimp} we deduce that $T_\nu^*=+\infty$. The wave front tracking is then globally defined in time.
\\
To summarize we have obtained uniform bound in $\nu$ on $z_\nu$ in $L^\infty_t(L^\infty)$ and on $w_\nu$ in $L^\infty_t(BV^{\frac{1}{3}}(\R))$, we wish now to develop some compactness argument in order to pass to the limit when $\nu$ goes to $+\infty$. The difficulty is to prove in particular that $z_\nu$ converges strongly to $z$ in $L^1_{loc,t,x}$ since we can not use Helly Theorem as it is the case for $w_\nu$.
\subsubsection*{Compactness argument for $(z_\nu)_{\nu>0}$}
We consider now the Lipschitz homeomorphism:
$$\phi^\nu (t,x)=(t,\gamma^\nu_2(t,x)),$$
with $\gamma^\nu_2(t,x)$ the forward generalized 2-characteristic such that $\gamma^\nu_2(0,x)=x$.
Furthermore we define $\eta^\nu_{L}$ and $z^\nu_L$ as follows:
$$\eta_\nu(t,x)=z_\nu(t,\gamma^\nu_2(t,x))-z_0(x)\;\;\mbox{and}\;\;z^\nu_L(t,x)=z_\nu(t,\gamma^\nu_2(t,x)).$$
We observe in particular that:
$$z_\nu(t,x)=z^\nu_L((\phi^{\nu})^{-1}(t,x)),$$
with $(\phi^{\nu})^{-1}$ the inverse  of the
 Lipschitz homeomorphism
$\phi^\nu$ (for instance see \cite{Giusti, Morgan}  for the notion of a Lipschitz homeomorphism and bi-Lipschitz homeomorphism, when the inverse is also Lipschitz, in geometric measure theory).
We are going now to prove a succession of different Lemmas.
\begin{lemme}
\label{lem3}
Up to a subsequence, we have:
$$ \lim_{\nu\rightarrow+\infty} z^\nu_L= z_L\;\;\mbox{in}\;L^1_{loc,t,x}.$$
\end{lemme}
{\bf Proof:} We have seen that $z_\nu$ is $BV$ along the curve $\{(t,\gamma^\nu_2(t,x)),\;t\in[0,T]\}$, indeed we have seen that  $z_\nu$ is increasing along $\{(t,\gamma^\nu_2(t,x)),\;t\in[0,T]\}$ and that $z_\nu$ is uniformly bounded in $L^\infty_{t,x}$ then $z_\nu$ is uniformly bounded in $\nu$ in $BV$ along the curve $\{(t,\gamma^\nu_2(t,x)),\;t\in[0,T]\}$. It implies that $z_L^\nu$ is uniformly bounded in $\nu$ in $L^\infty_x(BV([0,T]))$ for any $T>0$. Now we know that the speed of propagation of $z^\nu_L$ in the plane $(x,t)$ is finite, indeed the speed is bounded by:
\begin{equation}
\sup_{(t,x)\in S_{1}} \left |\frac{1}{\lambda_1((w^\nu(t,x),z^\nu(t,x))^-, (w^\nu(t,x),z^\nu(t,x))^+  )}\right |<M,
\label{vitesse}
\end{equation}
with $S_1$ the set of the 1-shock and $\lambda_1((w^\nu(t,x),z^\nu(t,x))^-, (w^\nu(t,x),z^\nu(t,x))^+  )$ the speed of the shock defined by the Rankine Hugoniot relation. Here
$M$ does not depend on $\nu$ and (\ref{vitesse}) is true because on a small square $[-r,r]^2$ with $r>0$ sufficiently small we have:
$$\sup_{x\in [-r,r]^2}\lambda_1(x)<0.$$
It is well known  \cite{Smoller, Lipconv} that it implies that $z^\nu_L$ is uniformly bounded in $Lip_x(L^1_{loc,t})$. Using the Kolmogorov theorem, we deduce that up to a subsequence in $\nu$  $z^\nu_L$ converges to $z$ in $L^1_{loc,t,x}$:
$$ \lim_{\nu\rightarrow+\infty} z^\nu_L = z_L\;\;\mbox{in}\;L^1_{loc,t,x}.$$
\cqfd
Next we use the following lemma.
\begin{lemme}
\label{diffa}
Let $\psi^\nu$ a homeomorphism uniformly Lipschitz in $\nu$ from $\R^+\times\R$ to $\R^+\times\R$ such that there exists $M>1$ verifying for any $\nu>0$ and {almost everywhere}:
\begin{equation}
0\leq\frac{1}{M}\leq |\det D_{t,x}\psi^\nu|\leq M,\;0\leq \| D_{t,x}\psi^\nu\|\leq M.
\label{conddiffeo}
\end{equation}
We assume that 
\begin{equation}
 \lim_{\nu\rightarrow+\infty} y^\nu = y \mbox{ in } L^1_{t,x,loc}, 
\mbox{ and }
 \lim_{\nu\rightarrow+\infty} \psi^\nu = \psi  \mbox{ in } L^\infty_{t,x,loc}
  \label{condibon}
  \end{equation}
   then,
\begin{align*} 
 \lim_{\nu\rightarrow+\infty} 
  y^\nu (\psi^\nu)=  y (\psi)  \mbox{ in } L^1_{t,x,loc}.
  \end{align*}
\end{lemme}

{\bf Proof:}  
Let $\va\in C^\infty_c(\R^+\times\R)$ a positive regular function with compact support and with values in $[0,1]$. For $\e>0$ we take $\widetilde{y}$ a continuous function in $L^1_{t,x}$ such that:
\begin{equation}
\|\widetilde{y}- y\|_{L^1_{t,x}(K_\va)}\leq\e.
\label{avi}
\end{equation}
with $K_\va$ a compact of $\R^+\times\R$ sufficiently large such that for any $\nu>0$ we have $\mbox{ supp }\va \left((\psi^\nu)^{-1}\right)$ which is included in $K_\va$.
 Using (\ref{conddiffeo}) and (\ref{avi}) we have:
 $$
\begin{aligned}
&\int_{\R^+}\int_{\R}\va |y^\nu \circ\psi^\nu-y \circ \psi| dx dt  \\ 
\leq  & \int_{\R^+}\int_{\R}\va |y^\nu \circ \psi^\nu-y \circ \psi^\nu | dx dt+\int_{\R^+}\int_{\R}\va|y  \circ \psi^\nu-\widetilde{y} \circ\psi^\nu | dx dt\\
&+\int_{\R^+}\int_{\R} \va |\widetilde{y} \circ\psi^\nu-\widetilde{y} \circ\psi | dx dt
+\int_{\R^+}\int_{\R}\va |\widetilde{y}  \circ \psi-y \circ \psi | dx dt\\
\leq  & 2\e M+M\|y-y^{\nu}\|_{L^1_{t,x}(K_\va)}+\int_{\R^+}\int_{\R}\va |\widetilde{y} \circ \psi^\nu-\widetilde{y} \circ \psi | dx dt.
\end{aligned}
$$
Using dominated convergence we can deal with the last integral and prove that for any positive $\va\in C^\infty_c(\R^+\times\R)$ , we have  $\va y^\nu \circ \psi^\nu\rightarrow \va y  \circ \psi$ in $L^1_{t,x}$. We deduce then that  $ y^\nu \circ \psi^\nu\rightarrow  y \circ \psi$ in $L^1_{t,x,loc}$
\cqfd
\begin{lemme}
\label{lem5}
$\phi^\nu$  and $(\phi^\nu)^{-1}$
verify the assumption (\ref{conddiffeo}) of the Lemma \ref{diffa}.
\end{lemme}
{\bf Proof:} It suffices to verify that there exist $M>1$ such that:
\begin{equation}
0\leq\frac{1}{M}\leq |\det D_{t,x}\phi^\nu|\leq M,\;0\leq \| D_{t,x}\phi^\nu\|\leq M.
\label{phi}
\end{equation}
We observe that:
\begin{equation}
D_{t,x}\phi^\nu(t,x)=\quad
\begin{pmatrix} 
1 & 0 \\
\p_t \gamma^\nu_2(t,x) & \p_x\gamma^\nu_2(t,x)
\end{pmatrix}
\label{det1}
\end{equation}
and:
\begin{equation}
D_{t,x}(\phi^\nu)^{-1}(\phi^\nu(t,x))=\quad
\begin{pmatrix} 
1 & 0 \\
\displaystyle
-\frac{\p_t \gamma^\nu_2(t,x)}{ \p_x\gamma^\nu_2(t,x)} 
& \displaystyle \frac{1}{ \p_x\gamma^\nu_2(t,x)}
\end{pmatrix}
\label{det2}
\end{equation}
We know that:
$$\p_t\gamma^\nu_2(t,x)\in[\min_\pm(\lambda_2(w^\nu(t,\gamma^\nu_2(t,x)^\pm))),\max_\pm(\lambda_2(w^\nu(t,\gamma^\nu_2(t,x)^\pm)))]$$
It implies in particular that $\p_t\gamma^\nu_2$ is uniformly bounded in $\nu$ since we have seen that $w_\nu$ is uniformly bounded in $L^\infty_{t,x}$. Similarly for  $\p_x\gamma^\nu_2(t,x)$, we observe that $v_\nu(t,x)=\p_x\gamma^\nu_2(t,x)$ verifies:
$$\p_t v_\nu(t,x)=\p_x(\lambda_2(w_\nu))(t,\gamma_2^\nu(t,x)) v_\nu(t,x) ,\;v_\nu(0,x)=x.$$
Indeed we can observe in fact that $\gamma^\nu_2$ verify except at the point where $\gamma_2^\nu$ meets a 1-wave front:
$$\p_t \gamma^\nu_2(t,x)=\lambda_2(w_\nu(t,\gamma_2^\nu(t,x))).$$
We deduce that:
\begin{equation}
v_\nu(t,x)=\exp\left(\int^t_0 \p_x(\lambda_2(w_\nu))(s,\gamma_2^\nu (s,x))ds \right).
\label{defiu}
\end{equation}
and:
\begin{equation}
\det D_{t,x}\phi^\nu=v_\nu(t,x).
\label{detiu}
\end{equation}
We have now:
\begin{equation}
\int^t_0 \p_x(\lambda_2(w_\nu))(s,\gamma_2^\nu(s,x))ds=\sum_{\alpha\in J_1}[\lambda_2(w_{\alpha,\nu})].
\label{teles}
\end{equation}
$J_1$ is the set of point where a 1-wave front meets the curve $\{(\alpha t,\gamma^\nu_{2,\alpha}(t)),\;\;t\in[0,T]\}$. Since in (\ref{teles}), we have a telescopic sum, we deduce that:
\begin{equation}
-2\|\lambda_2(w_\nu(\cdot,\cdot))\|_{L^\infty_{t,x}}\leq \int^t_0 \p_x(\lambda_2(w_\nu))(s,\gamma_2^\nu(s,x))ds\leq 2\|\lambda_2(w_\nu(\cdot,\cdot))\|_{L^\infty_{t,x}}
\label{teles1}
\end{equation}
From (\ref{defiu}) and (\ref{teles1}) we deduce that $\p_x\gamma_2^\nu$ and $\frac{1}{\p_x\gamma_2^\nu}$ is uniformly bounded in $\nu$. Furthermore it implies also that the determinant of $D_{x,t}\phi^\nu$ and $D_{x,t}(\phi^\nu)^{-1}$ satisfies the assumption  (\ref{phi}) using (\ref{detiu}). Since $\p_t\gamma_2^\nu$ is uniformly bounded in $\nu$ we deduce finally that 
$\phi^\nu$ and $(\phi^\nu)^{-1}$verifies uniformly in $\nu$ (\ref{phi}) using the formula (\ref{det1}) and (\ref{det2}).
\cqfd
\begin{lemme}
\label{lem6}
$((\phi^\nu)^{-1})_{\nu>0}$ converges up to a subsequence to $\phi^{-1}$ in $L^\infty_{t,x,loc}$.
\end{lemme}
{\bf Proof:} Indeed for any compact $K$ of $\R^+\times\R$ we observe that $(\phi^\nu)^{-1}$ is a continuous function from $K$ to $\R^+\times\R$ because $(\phi^\nu)^{-1}$ is Lipschitz using the Lemma \ref{lem5}. Since the sequence $((\phi^\nu)^{-1})_{\nu>0}$ is uniformly Lipschitz in $C(K,\R^+\times\R)$ from the Lemma \ref{lem5}, using the Ascoli Theorem up to a subsequence $((\phi^\nu)^{-1})_{\nu>0}$ converges uniformly to $\phi^{-1}$ on $K$. Using a standard argument of diagonal process we obtain that $((\phi^\nu)^{-1})_{\nu>0}$ converges up to a subsequence to $(\phi)^{-1}$ in $L^\infty_{t,x,loc}$.
\cqfd
Using the Lemmas \ref{lem3}, \ref{diffa}, \ref{lem5} and  \ref{lem6} we deduce that $z^\nu_Lo(\phi^\nu)^{-1}=z_\nu$ converges strongly to $z=z_L(\phi^{-1})$  in $L^1_{loc,t,x}$ up to a subsequence. It implies that up to a subsequence $z_\nu$ converges almost everywhere up to a subsequence to $z=z_L(\phi^{-1})$ when $\nu$ goes to $+\infty$. concerning the convergence of the sequel  $(w_\nu)_{\nu>0}$, the proof is more simple since  $(w_\nu)_{\nu>0}$ is uniformly bounded in $L^\infty(\R^+,BV^{\frac{1}{3}})$. 
Thus, as in \cite{BGJ6} for the scalar case or \cite{BGJP} for a $2\times 2$ system $w^\nu$ is also bounded in $Lip^s_t([0,+\infty[,L^p_{loc}(\R,\R))$
with $p=1/s$ and the compactness follows. \\
We deduce then that the sequence $u_{\nu}=(w,z)^{-1}(w_{\e_\nu},z_{\e_\nu})$ with $(w,z)^{-1}$ the inverse of the local diffeomorphism $(w,z)$ converges also almost everywhere to $u$ with $u_\nu$ uniformly bounded in $L^\infty_{t,x}$. It is then classical to verify that $u$ is a global weak solution of the system (\ref{1}) using dominated convergence (see \cite{Br,Da,Smoller} for more details).

\subsubsection*{Decomposition of $z$ in Theorem \ref{thm:1/3,0}}
 Finally, the decomposition of $z$  \eqref{zdecomposition} for the weak solution can be proved.    We already use  this decomposition for the approximate sequence $(z^\nu)_{\nu>0}$ with: 
$$
z_\nu(t,\gamma_2^\nu(t,x))= z_{0,\nu}(x) + \eta_\nu(t,x).
$$ 
 The 2-characteristics $(\gamma_2^\nu)$  are equi-Lipschitz. Thus, up to  a subsequence, we can pass to the limit when $\nu \rightarrow + \infty$. At the limit, $\gamma_2$ satisfies  the differential equation for the generalized 2-characteristics \cite{Da,Fi}. Moreover, $z_{0,\nu}$ and $\eta_\nu$ converge in $L^1_{loc}$, the whole sequence for $(z_{0,\nu})$ and only a subsequence for $(\eta_\nu)$.  Now, we can pass to the limit in $ z_{\nu}(t,x) $ using the bi-Lipschitz homeomorphism, as previously with Lemma \ref{diffa}, 
 to obtain the decomposition  \eqref{zdecomposition}.
The proof of Theorem \eqref{thm:1/3,0} is achieved.
\section*{Acknowledgements}
Boris Haspot has been partially funded by the ANR project INFAMIE ANR-15-CE40-0011. This
work was partially realized during the secondment of Boris Haspot in the ANGE Inria team.
 

\end{document}